\documentclass[11pt]{article}
\usepackage[centertags]{amsmath}
\usepackage{amsfonts}
\usepackage{amssymb}
\usepackage{amsthm}
\usepackage{amsmath}
\usepackage{version}
\usepackage{tabularx}
\usepackage{bbm}
\usepackage{eurosym}
\usepackage{placeins}
\usepackage{subfigure}
\usepackage{slashbox}
\usepackage[german,english]{babel}
\usepackage{latexsym,graphicx}
\usepackage{color}
\usepackage{enumitem}

\theoremstyle{plain}
\newtheorem{thm}{Theorem}[section]

\newtheorem{prop}[thm]{Proposition}
\newtheorem{assump}{Assumption}

\theoremstyle{definition}
\newtheorem{defn}{Definition}[section]

\theoremstyle{remark}
\newtheorem{rem}{\bf Remark}[section]
\theoremstyle{remark}

\usepackage{fancyhdr}

\DeclareMathOperator*{\Argmin}{Arg\,min}

\makeatletter
\def \newequation#1#2{
   \@definecounter{#1}
   \@namedef{the#1}{\hbox{#2}}
   \@namedef{#1}{$$\refstepcounter{#1}}
   \@namedef{end#1}{
      \eqno \csname the#1\endcsname $$\global\@ignoretrue
      }
}
\makeatother
\newequation{E1}{($E_{b,\sigma,W}$)}
   
\makeatletter
\def \newequation#1#2{
   \@definecounter{#1}
   \@namedef{the#1}{\hbox{#2}}
   \@namedef{#1}{$$\refstepcounter{#1}}
   \@namedef{end#1}{
      \eqno \csname the#1\endcsname $$\global\@ignoretrue
      }
   }
\makeatother
\newequation{hyp3}{($\mathcal{H}_{b,\sigma}$)}

\makeatletter
\def \newequation#1#2{
   \@definecounter{#1}
   \@namedef{the#1}{\hbox{#2}}
   \@namedef{#1}{$$\refstepcounter{#1}}
   \@namedef{end#1}{
      \eqno \csname the#1\endcsname $$\global\@ignoretrue
      }
   }
\makeatother
\newequation{shleg}{(ShLeg)}

\makeatletter
\def \newequation#1#2{
   \@definecounter{#1}
   \@namedef{the#1}{\hbox{#2}}
   \@namedef{#1}{$$\refstepcounter{#1}}
   \@namedef{end#1}{
      \eqno \csname the#1\endcsname $$\global\@ignoretrue
      }
   }
\makeatother
\newequation{kl}{(KL)}

\makeatletter
\def \newequation#1#2{
   \@definecounter{#1}
   \@namedef{the#1}{\hbox{#2}}
   \@namedef{#1}{$$\refstepcounter{#1}}
   \@namedef{end#1}{
      \eqno \csname the#1\endcsname $$\global\@ignoretrue
      }
   }
\makeatother
\newequation{haar}{(Haar)}

\begin{document}

\vspace{\stretch{1}}

\title{Quantization based Recursive Importance Sampling}
          
\author{ 
{\sc  Noufel Frikha} \thanks{CMAP, Ecole Polytechnique Paris, CNRS, 
e-mail: {\tt noufel.frikha@cmap.polytechnique.fr} } \ \ \ 
{\sc  Abass Sagna} \thanks{Laboratoire d'Analyse et de Probabilit\'es de l'Universit\'e d'Evry Val d'Essonne \& ENSIIE. The author's research is supported by an AMaMeF exchange grant and the ``Chaire Risque de Cr\'{e}dit'' of the French Banking Federation. 
e-mail: {\tt abass.sagna@gmail.com} } 
}
\maketitle
\begin{abstract}
We investigate in this paper an alternative method to simulation based recursive importance sampling procedure to estimate the optimal change of measure for Monte Carlo simulations. We propose an algorithm which combines (vector and functional) optimal quantization with Newton-Raphson zero search procedure. Our approach can be seen as a robust and automatic deterministic counterpart of recursive importance sampling by means of stochastic approximation algorithm which, in practice, may require tuning and a good knowledge of the payoff function in practice. Moreover, unlike recursive importance sampling procedures, the proposed methodology does not rely on simulations so it is quite generic and can come along on the top of Monte Carlo simulations.

We first emphasize on the consistency of quantization for designing an importance sampling algorithm for both multi-dimensional distributions and diffusion processes. We show that the induced error on the optimal change of measure is controlled by the mean quantization error. 

We illustrate the effectiveness of our algorithm by pricing several options in a multi-dimensional and infinite dimensional framework.
  
\end{abstract}
   

%
\section*{Introduction} 
In this paper, we are interested in one of the most basic problem of numerical probability which consists in the computation of the expectation 
\begin{equation}
\mathbb{E}[F(X)]
\label{eq0}
\end{equation}

\noindent where $X:(\Omega, \mathcal{A}, \mathbb{P}) \rightarrow \left(E, |.|_{E}\right)$ is a random vector taking values in a Banach space $E$ and $F: E \rightarrow \mathbb{R}$ is a Borel function such that $\mathbb{E}[F(X)^{2}]<+\infty$. When the space $E$ is $\mathbb{R}^{d}$ (equipped with the Euclidean norm) we will refer to the finite (multi) dimensional setting and when the space $E$ is $\mathcal{C}([0,T],\mathbb{R}^{d})$ (equipped with the supremum norm) to deal with the case where $X$ is a continuous path-dependent diffusion process, we will refer to the infinite dimensional setting. For instance in mathematical finance, computing the price of an option and the sensitivities of this price with respect to some parameters amounts to estimate such a quantity. When no closed or semi-closed formulas are available, one often relies on Monte Carlo simulation which remains the most widely used numerical method in this context. 

Variance reduction methods are often used to increase the accuracy of a Monte Carlo simulation or reduce its computation time. The most common variance reduction methods are antithetic variables, conditioning, control variate, importance sampling and stratified sampling. Adaptive variance reduction methods have been recently investigated to take advantage of the random samples used to compute the expectation above in order to optimize at the same time the variance reduction tool (see \cite{Jourdain2009a} and the references therein). In practice, it is not clear that this adaptive one step Monte Carlo procedure is better than the basic two step procedure: optimizing the variance reduction tool (e.g. the optimal change of measure in an importance sampling framework) using a small number of samples, then computing the expectation of interest with this optimized parameter. 

In this paper, we are interested in variance reduction by importance sampling (IS). We denote by $p$ the probability density function of $X$ and we consider a parameterized family of density functions $\left(p_{\theta}\right)_{\theta \in \Theta}$: we will choose $\Theta=\mathbb{R}^{d}$ in the finite dimensional setting and $\Theta = L^{2}([0,T],\mathbb{R}^{q})$, where $q$ is the dimension of the Brownian motion driving the dynamic of the $d$-dimensional process $X$, when dealing with the infinite dimensional setting. This general framework is the one investigated in \cite{Lemaire2008}. In this introduction we will present the importance sampling paradigm in the finite multi-dimensional case, $i.e.$ we set $E=\mathbb{R}^{d}$. We suppose that $p>0$, $\lambda_{d}-a.e.$, where $\lambda_{d}$ denotes the Lebesgue measure on $(\mathbb{R}^{d},\mathcal{B}or(\mathbb{R}^{d}))$ and we set $p_{0}=p$. Moreover, we focus throughout all the paper on importance sampling by mean translation, $i.e.$ we will consider that for all $\theta \in \mathbb{R}^{d}$, for all $x \in \mathbb{R}^{d}$, $p_{\theta}(x)=p(x-\theta)$. The basic idea of IS is to introduce the parameter $\theta$ in the above expectation \eqref{eq0} using the invariance by translation of the Lebesgue measure, for every $\theta \in \mathbb{R}^{d}$
\begin{equation}
\mathbb{E}[F(X)] = \mathbb{E}\left[F(X+\theta)\frac{p(X+\theta)}{p(X)}\right],
\label{eq1}
\end{equation}     

\noindent and among all these random variables with the same expectation, we want to select the one with the lowest variance, $i.e.$ the one with the lowest quadratic norm
$$
Q(\theta):=\mathbb{E}\left[F^{2}(X+\theta)\frac{p^{2}(X+\theta)}{p^{2}(X)}\right], \ \ \theta \in \mathbb{R}^{d}.
$$

\noindent A reverse change of variable shows that:
\begin{eqnarray}
Q(\theta)=\mathbb{E}\left[F^{2}(X) \frac{p(X)}{p(X-\theta)}\right],  \hspace*{0.2cm} \theta \in \mathbb{R}^{d}.
\label{eq2}
\end{eqnarray}

\noindent Now if the density function $p$ of $X$ satisfies
\begin{equation}
(i) \hspace*{0.2cm} p \mbox{ is log-concave } \ \ \mbox{ and } \ \
                            (ii) \hspace*{0.2cm} \lim_{|x| \rightarrow +\infty } p(x)=0                       
\label{hyp1}
\end{equation} 

\noindent and
\begin{equation}
 Q(\theta) < +\infty, \mbox{ } \forall \theta \in \mathbb{R}^{d}                
\label{hyp2}
\end{equation} 

\noindent then, one shows that the function $Q$ is finite, convex and goes to infinity at infinity, thus $\arg \min_{\theta} Q={\left\{\theta \in \mathbb{R}^{d} \ | \ D Q(\theta) =0\right\}}$, where $D Q$ is the gradient of $Q$, is non empty (for a proof, we refer to \cite{Lemaire2008}). 

\noindent Now, if $D Q$ admits a representation as an expectation, then it is possible to devise a recursive Robbins-Monro (RM) procedure to approximate the optimal parameter $\theta^{*}$ minimizing $Q$, namely
\begin{equation}{}
\theta_{n}=\theta_{n-1}-\gamma_{n}K(\theta_{n-1},X_{n}), \ n\geq1
\label{eq3}
\end{equation}

\noindent where $(X_{n})_{n\geq1}$ is an $i.i.d.$ sequence of random vectors having the distribution of $X$, $(\gamma_{n})_{n\geq1}$ is a positive deterministic sequence satisfying,
$$
\sum_{n\geq1} \gamma_{n} = +\infty, \ \ \mbox{ and } \sum_{n\geq1} \gamma_{n}^{2} < +\infty, 
$$
 
\noindent and $K$ is naturally defined by the formal differentiation of $Q$: for every $x \in \mathbb{R}^{d}$,
\begin{equation}
K(\theta,x) = F^{2}(x)\frac{p(x)}{p^{2}(x-\theta)}D
p(x-\theta).
\label{eq4}
\end{equation}

\noindent IS by means of stochastic approximation has been investigated by several authors (see e.g. \cite{Glynn1989}, \cite{Dufresne1998} and \cite{Fu2000}) in order to estimate the optimal change of measure by a RM procedure. It has recently been studied in the Gaussian framework in \cite{Arouna2004} (see also \cite{ArounaWinter2003/04}) where \eqref{eq4} is used to design a stochastic gradient algorithm. However, the regular RM procedure \eqref{eq4} suffers from an instability issue coming from the fact that the classical sub-linear growth Assumption in quadratic mean in the Robbins-Monro Theorem 
\begin{equation}
\forall \theta \in \mathbb{R}^{d}, \ \ \mathbb{E}\left[K(\theta,X)^{2}\right]^{\frac{1}{2}} \leq C\left(1+|\theta|\right)
\label{eq7}
\end{equation}

\noindent is only fulfilled when $F$ is constant, due to the behaviour of the annoying term $p(x)/p(x-\theta)$ as $\theta$ goes to infinity. Consequently, $\theta_{n}$ can escape at infinity at almost every implementation as pointed out in \cite{Arouna2004}. To circumvent this problem, a ``projected version'' of the procedure based on repeated reinitializations when the algorithms exits from an increasing sequence of compact sets (while the step $\gamma_{n}$ keeps going to 0) was used. This approach is known as the projection ``à la Chen'' (see for instance \cite{Chen1987}, \cite{Chen1986} and \cite{Lelong2008}). A central limit theorem for this version of the recursive importance sampling algorithm is proved in \cite{Lelong2007}. This kind of technique forces the stability of the algorithm and prevents explosion. From a numerical point of view, this projected algorithm is known to converge if the sequence of compact sets have been specified suitably, which is not an easy task in practice. 

In \cite{Lemaire2008}, the authors propose a third change of variable to plug back the parameter $\theta$ into $F$ 
\begin{equation}
D Q(\theta) = \mathbb{E}\left[K(\theta,X)\right] = \mathbb{E}\left[F^{2}(X-\theta)\frac{p^2(X-\theta)}{p(X)p(X-\theta)}\frac{D
p(X-2\theta)}{p(X-2\theta)}\right],
\label{eq8}
\end{equation}

\noindent which has a known controlled growth rate at infinity in common applications. For instance, if we suppose that there exists $\lambda>0$, such that 
$$
\forall \ x \in \mathbb{R}^{d}, \ \ \left|F(x)\right| \leq C \ e^{\lambda|x|},
$$

\noindent we can define a new function $H$ by setting
\begin{equation}
H(\theta,x) := e^{-2\lambda |\theta|} F^{2}(x-\theta)\frac{p^2(x-\theta)}{p(x)p(x-\theta)}\frac{D
p(x-2\theta)}{p(x-2\theta)}
\label{eq9}
\end{equation}

\noindent Under additional assumptions on the density function $p$, they derive a new stochastic approximation algorithm using $H$ which satisfies the sub-linear growth condition \eqref{eq7} so that it $a.s.$ converges and is stable without having to apply projection technique. They also extend this construction to exponential change of measure (Esscher transform) and to diffusion process using Girsanov transform. However, from a numerical point of view, the tuning of the algorithm needs a good knowledge of the behavior of $F$ at infinity. In practical implementations, recursive importance sampling methods using stochastic implies specific tuning of the step sequence or the sequence of compact sets which in both cases strongly depends of the payoff function $F$.
  
In order to get rid of this problem, \cite{Jourdain2009} proposes an optimization Newton's algorithm to estimate the optimal change of measure in a Gaussian framework. They approximate $D Q(\theta)$ and $D^{2} Q(\theta)$ (based on the representation \eqref{eq4} with $p(x)=\frac{1}{\sqrt{2\pi}}e^{-\frac{|x|^{2}}{2}}$, $x\in \mathbb{R}^{d}$) using Monte Carlo simulation with $n$ samples
\begin{align*}
D Q_{n}(\theta)  & = \frac{1}{n} \sum_{k=1}^{n} \left(\theta-X_{k}\right)F^{2}(X_{k})e^{-\theta.X_{k}+\frac{|\theta|^{2}}{2}}   \\
D^{2} Q_{n}(\theta)   & = \frac{1}{n} \sum_{k=1}^{n} \left(I_{d}+\left(\theta-X_{k}\right)\left(\theta-X_{k}\right)^{T}\right)F^{2}(X_{k})e^{-\theta.X_{k}+\frac{|\theta|^{2}}{2}}
\end{align*}
 
\noindent where $(X_{k})_{1\leq k \leq n}$ is an $i.i.d.$ sequence of $d$-dimensional standard normal vectors. The (unique) minimum of $Q$ is approximated by the (unique) zero $\theta_{n}$ of $Q_{n}$ which can be computed using a deterministic Newton-Raphson algorithm. Moreover, several asymptotic properties are addressed.

Deterministic optimization using a large deviation argument has been investigated in \cite{Glasserman1999}. The optimal change of measure is selected as a local maximum of $\theta \mapsto \log |F(\theta)| - \frac{|\theta|^{2}}{2}$. This can be achieved only under some regularity assumptions on the function $F$. Moreover, from a theoretical point of view this choice is not optimal.

In this work, we propose and study an alternative deterministic procedure which is not adaptive and does not need specific payoff-dependent tuning. The optimal change of measure is estimated using a robust and automatic Newton-Raphson's algorithm combined with optimal (vector and functional) quantization. The main advantage of our procedure is that it can be used as a generic variance reduction method which comes upstream the Monte Carlo simulation framework. Moreover, we will focus on one very common situation in mathematical finance, that is, Monte Carlo simulation that are based on multi-factor Brownian diffusions. Indeed, the presented method easily extends to this framework using quadratic optimal functional quantization of stochastic processes (see $e.g.$ \cite{Pag`es2007}). It is particularly adapted for financial institutions since the methodology we propose can come along on the top of Monte Carlo simulations. Numerical tests ilustrate the effectiveness of our approach in both multi-dimensional and infinite dimensional frameworks. 

The paper is organized as follows. Section 1 presents several results about vector and functional quantization that are required in the following. We focus on the functional quantization of diffusion processes. Section 2 presents the quantization based recursive importance sampling algorithm. The emphasis is on the consistency of quantization for designing an importance sampling algorithm for both multi-dimensional distributions and diffusion processes. We show in particular that the induced error on the optimal change of measure is controlled by the mean quantization error. In section 3, we provide numerical experiments of our approach by considering option pricing problems arising in mathematical finance.  

\section{Some results on optimal quantization}
Before, dealing with the construction of the quantization based IS algorithm, we provide with some background on quantization of Hilbert spaces and Gaussian processes viewed as $L^{2}_{T}$-valued random vectors. 

\subsection{Introduction to quantization of random variables}
Let $N \in \mathbb{N}^{*}$. The principle of the $N$-quantization of a random variable $X$ taking its values in a separable Hilbert space $E$ is to study the best $\left\|.\right\|_{p}$-approximation of $X$ by $E$-valued random vectors taking at most $N$ values. The norm $\left\|.\right\|_{p}$ is the usual norm on $L^{p}_{E}(\Omega,\mathbb{P})$ defined by $\left\|X\right\|_{p}=\left(\mathbb{E}\left[|X|_{E}^{p}\right]\right)^{1/p}$. When $p=2$, we talk about quadratic optimal quantization. If $E=\mathbb{R}^{d}$, one speaks about vector quantization. When $E$ is an infinite dimensional Hilbert space like $L^{2}\left([0,T],dt\right)$ endowed with the usual norm $|f|_{E}=\left(\int_{0}^{T} f(t)^{2} dt\right)^{\frac{1}{2}}$, we talk about functional quantization. 

\begin{defn}[Voronoi tessellation]
Let $N\geq1$ and $x:=(x_{1},\cdots,x_{N}) \in E^{N}$ be a $N$-tuple referred to us by a $N-$quantizer and let Proj$_{x}: E \rightarrow \left\{x_{1},\cdots,x_{N}\right\}$ be a projection following the closest neighbour rule. Then, the Borel partition $C=\left\{C_{1},\cdots,C_{N}\right\}$ of $E$ defined by $C_{i}:=$Proj$^{-1}_{x}({x_{i}})$, $i=1,\cdots,N$ and satisfying
$$
\mbox{Proj}_{x}^{-1}(\left\{x_{i}\right\}) \subset \left\{ y \in E, |x_{i}-y|_{E} = \min_{1 \leq j \leq N}|x_{j}-y|_{E}\right\}, \ \ 1\leq i \leq N.
$$

\noindent is called a Voronoi tessellation of $E$ induced by $x$.
\end{defn}

\noindent One defines the Voronoi quantization of $X$ induced by $x$ by
$$
\widehat{X}^{x} = \mbox{Proj}_{x}\left(X\right)=\sum_{i=1}^{N} x_{i} \mbox{\textbf{1}}_{\left\{X\in C_{i}\right\}}.
$$

\noindent The discrete random variable $\widehat{X}^{x}$ (one will sometimes simply write $\widehat{X}$ if there is no ambiguity) is the best $L^{p}(\mathbb{P})$-approximation of $X$ among all measurable random variable taking values in $x:=\left\{x_{1},\cdots,x_{N}\right\}$. In fact, for any random variable $Y: \Omega \rightarrow \left\{x_{1}, \cdots, x_{N}\right\}$, we have 
$$
\left|X-\widehat{X}\right|_{E} \leq \left|X-Y\right|_{E} \ \ \mathbb{P}-a.s. 
$$

\noindent so that $\left\|X-\widehat{X}\right\|_{p} \leq \left\|X-Y\right\|_{p}$.

For any fixed $N$-quantizer $x:=(x_{1},\cdots,x_{N})$, we associate the $L^{p}(\mathbb{P})$-mean error $\left\|X-\widehat{X}\right\|_{p}$ induced by $x$. One aims at finding a $N-$tuple $x\in E^{N}$ which minimizes the $L^{p}$-mean error over $E^{N}$. It amounts to minimizing the function
$$
Q_{N}^{X}:\left(x_{1},\cdots,x_{N}\right) \mapsto \left\|X-\hat{X}\right\|_{p}=\left\|\min_{1 \leq i \leq N }   |X-x_{i}|_{E}\right\|_{p}
$$

The function $Q_{N}^{X}$ reaches its minimum at one (at least) tuple $x^{*}$ called an optimal $N$-quantizer. This infimum is in general not unique, except in some cases, in particular when $d=1$ and the density of $X$ is log-concave. Note that card($x^{*}$)$=N$ if card(supp($\mathbb{P}_{X}$))$\geq N$. Moreover, the $L^{p}$-mean quantization error $e_{N,p}^{X}:=\min Q_{N}^{X}$ converges toward $0$ and for non-singular $\mathbb{R}^{d}$-valued random vectors, the rate of convergence of convergence of $e_{N,p}^{X}$ is ruled by the so-called Zador Theorem (see \cite{Graf2000}).
\begin{thm}
Assume that $X\in L^{p+\delta}_{\mathbb{R}^{d}}(\mathbb{P})$ for some $\delta>0$. Let $f$ denote the density of the absolutely continuous part of $\mathbb{P}_{X}$ (possibly, $f\equiv0$). Then,
\begin{equation}
\lim_{N\rightarrow +\infty}N^{1/d}e_{N,p}^{X} = q_{p}(d) \left(\int f(x)^{d/d+p} dx\right)^{(d+p)/dp},
\label{eq10}
\end{equation}

\noindent where $q_{p}(d)$ is a strictly positive depending only on $p$ and on the dimension $d$.
\end{thm}

For further theoretical results on optimal vector quantization we refer to \cite{Graf2000}. One of the important issues for the Numerical Probabilities viewpoint is to compute the optimal quantizers and the associated weights (for the applications to Numerical Probabilities in the finite dimensional case, see the seminal paper \cite{Pag`es1998}). However, due to the non-uniqueness of the optimal quantizers in the general framework, specially when $E=\mathbb{R}^{d}$ with $d\geq2$, we are usually leaded to search for stationary quantizers, $i.e.$ $N$-quantizers $x$ satisfying $\nabla Q_{N}^{X}(x)=0$.  We will see further on that stationary quantizers are an important class of quantizers for numerics. The commonly used result is the quadratic case (when $p=2$) recalled below.

\begin{defn}[Stationarity]
Let $x:=(x_{1},\cdots,x_{N})\in E^{N}$ be an $N$-quantizer and $C:=\left\{C_{1},\cdots,C_{N}\right\}$ its associated Voronoi partition. The random vector $\hat{X}$ is called a stationary $N$-quantization of $X$ if it satisfies 
\begin{equation}
\forall \ i \neq j, \ \ x_{i}\neq x_{j} \ \ \mbox{and} \ \ \mathbb{P}\left(X \in \cup_{i} \partial C_{i}(x) \right)=0
\label{eq11}
\end{equation}

\noindent ($\mathbb{P}_{X}$-negligible boundary of the Voronoi cells) and
$$
\mathbb{E}\left[X \ | \ \widehat{X}\right] = \hat{X}.
$$

\noindent We have in particular, $\mathbb{E}[X] = \mathbb{E}[\hat{X}]$.
\end{defn}

\subsection{Some backgrounds on functional quantization}
A rigorous extension of optimal vector quantization to functional quantization is done in \cite{Luschgy2002}. The vector quantization problem is transposed to random vectors in an infinite dimensional Hilbert space, in particular, to stochastic processes $\left(X_{t}\right)_{t\in[0,T]}$ viewed as random vectors with values in $L^{2}\left([0,T],dt\right)$. In \cite{Pag`es2005}, numerical performances of quadratic functional quantization with applications to finance is investigated. In particular, the roles played by product quantizers and the so-called Karhunen-Lo\`eve (K-L) expansion of Gaussian processes are pointed out.

In what follows, we will start by the functional quantization of the standard Brownian motion since everything can be made explicit for this process. Then we will show how to construct from optimal quadratic functional quantization of Brownian motion explicit (non-Voronoi) quantization of Brownian diffusions. 

Assume that the separable Hilbert space $E$ is $L^{2}_{T}:=L^{2}([0,T],dt)$, with $\left\langle f,g \right\rangle=\int_{0}^{T}f(t)g(t) dt$. One defines the covariance operator $C_{W}$ of the Brownian motion $(W_{t})_{t\in[0,T]}$, for every $f\in L^{2}_{T}$ by
$$
C_{W}(f):=\mathbb{E}\left[\left\langle f,W\right\rangle W\right]=\left(t\mapsto\int_{0}^{T}f(s) (s \wedge t)ds\right).
$$

\noindent This operator is symmetric positive and can be diagonalized in the K-L orthonormal basis of $L^{2}_{T}$ $(e_{n})_{n\geq1}$ with eigenvalues $(\lambda_{n})_{n\geq1}$ given by
$$
e_{n}(t)=\sqrt{\frac{2}{T}} \sin\left(\pi \left(n-\frac{1}{2}\right)\frac{t}{T}\right), \ \ \lambda_{n}=\left(\frac{T}{\pi\left(n-\frac{1}{2}\right)}\right)^{2}, \ \ n \geq 1.
$$

\noindent Moreover, one may expand the paths of $(W)_{t \in [0,T]}$ on this basis, $i.e.$
\begin{equation}
W \stackrel{L^{2}_{T}}{=} \sum_{n \geq 1} \left\langle W, e_{n}\right\rangle e_{n}, \ \ \ \mathbb{P}\mbox{-}a.s. 
\label{eq12}
\end{equation}

\noindent Using Fubini's Theorem and the orthonormality of the K-L basis, one obtains for $\ell \geq 1$, $p \geq 1$
$$
\mathbb{E}\left[\left\langle W,e_{\ell}\right\rangle \left\langle W,e_{p}\right\rangle\right] = \left\langle C_{W}(e_{\ell}),e_{p}\right\rangle=\lambda_{\ell} \delta_{\ell p}
$$

\noindent where $\delta_{\ell p}$ denotes the Kronecker symbol. Consequently, the Gaussian sequence $\left(\left\langle W,e_{\ell}\right\rangle\right)_{\ell \geq 1}$ is pairwise non-correlated so that these random variables are independent. Hence, \eqref{eq12} can be written
$$
W \stackrel{L^{2}_{T}}{=} \sum_{n \geq 1} \sqrt{\lambda_{k}} \xi_{n} e_{n},
$$

\noindent where $\xi_{n}:=\left\langle W,e_{n}\right\rangle /\sqrt{\lambda_{n}}$, $n\geq1$, is an i.i.d. sequence of random variables with standard normal distribution. 

Now the idea of product functional quantization using at most $N$ elementary quantizers is to quantize these random coordinates $\xi_{n}$, $i.e.$ for every $n \geq 1$, one considers an optimal $N_{n}$-quantization ($N_{n}\geq1$) of $\xi_{n}$, denoted $\hat{\xi}_{n}$ where $\hat{\xi}_{n}:=\mbox{Proj}_{x_{n}}\left(\xi_{n}\right)$, $x_{n}:=\left(x_{1}^{N_{n}},\cdots,x_{N_{n}}^{N_{n}}\right)$ is the unique optimal $N_{n}$-quantizer of the normal distribution and $N_{1} \times \cdots \times N_{n} \leq N$, $N_{1}, \cdots, N_{n}\geq1$. For $n$ large enough, we set $N_{n}=1$, $\hat{\xi}_{n}=0$ (which is the optimal 1-quantization) and we define the product quantizer by (the finite sum)
$$
\widehat{W}_{t} = \sum_{n\geq1} \sqrt{\lambda_{n}} \hat{\xi}_{n} e_{n}(t). 
$$

\noindent The product quantizer $\chi$ that produces the above Voronoi quantization $\widehat{W}$ is defined by
$$
\chi_{\underline{i}}(t) = \sum_{n\geq1} \sqrt{\lambda_{n}} x_{i_{n}} e_{n}(t), \ \ \ \underline{i}=(i_{1},\cdots,i_{n}, \cdots)\in  \prod_{n\geq1} \left\{1,\cdots,N_{n}\right\}.
$$

\noindent and for every multi-index $\underline{i} \in \prod_{n\geq1} \left\{1,\cdots,N_{n}\right\}$, the associated Voronoi cell of $\chi$ is
$$
C_{\underline{i}}(\chi) = \prod_{n\geq1} \sqrt{\lambda_{n}} C_{i_{n}}(x_{n}).
$$

\noindent Moreover, from the independence of the normal random variables $(\xi_{n})_{n\geq1}$ the weights $\mathbb{P}\left(\widehat{W} = \chi_{\underline{i}} \right)$ can be computed explicitly
$$
\mathbb{P}\left(\widehat{W} = \chi_{\underline{i}}\right) = \prod_{n \geq 1} \mathbb{P}\left(\xi_{n} \in C_{i_{n}}(x_{n})\right).
$$

\noindent For numerical purposes, one may be interested by the theoretical rate of convergence for the quantization error of the Brownian motion and the stationarity of K-L product quantizer.
\begin{prop}[stationarity, see \cite{Pag`es2005}]
The product quantizer of the Brownian motion defined above is a stationary quantizer, $i.e.$
$$
\mathbb{E}\left[W \ | \ \widehat{W}\right] = \widehat{W} .
$$

\end{prop}

\begin{prop}[convergence rate, see \cite{Luschgy2002}]
For every $N\geq1$, there exists an optimal product quantizer of size at most $N$, denoted $\widehat{W}$ of the Brownian motion defined as the solution to the minimization problem
\begin{equation}
\min \left\{ \left\|  W- \widehat{W} \right\|_{2}, \ N_{1},\cdots,N_{n} \geq 1, \ d_{N}:=N_{1}\times \cdots \times N_{n} \leq N, N \geq 1  \right\}
\label{minprobfunctquant}
\end{equation}

\noindent Furthermore, these optimal product quantizer induces a rate optimal sequence, $i.e.$ 
$$
\left\|W-\widehat{W}\right\|_{2} \leq C \frac{T}{(\log N)^{\frac{1}{2}}}.
$$

\noindent for some real constant $C>0$
\end{prop}

To conclude this section, we shortly describe a constructive way to quantize scalar brownian diffusions (for more details see $e.g.$ \cite{Pag`es2005}). The rate is $O\left(\left(\log N\right)^{-\frac{1}{2}}\right)$ like for the Brownian motion as soon as the diffusion coefficient is not too degenerate. Consider the homogeneous Brownian diffusion process:
$$
d X_{t} = b(X_{t}) dt + \sigma(X_{t}) dW_{t}, \ X_{0}=x_{0},
$$

\noindent where $b$ and $\sigma$ are continuous on $\mathbb{R}$ with at most linear growth ($i.e.$ $|b(x)|+|\sigma(x)| \leq C(1+|x|)$) so that a weak solution to the equation exists. Let $(\chi^{N})_{N\geq1}$ be a sequence of rate-optimal K-L product quantizers of the Brownian motion. For every multi-index $\underline{i} \in \prod_{n\geq1} \left\{1,\cdots,N_{n}\right\}$, with $N_{1}\times \cdots \times N_{n} \leq N$, consider $x_{\underline{i}}$ the solution of the following integral equations
\begin{equation}
dx_{\underline{i}}(t) = \left(b(x_{\underline{i}}(t))-\frac{1}{2} \sigma \sigma'(x_{\underline{i}}(t))\right) dt + \sigma(x_{\underline{i}}(t)) d\chi_{\underline{i}}(t),
\label{eq13}
\end{equation}

\noindent where $\sigma'$ is the first derivative of $\sigma$. To simplify notations we consider the simpler notation $i$ for the multi-index $\underline{i}$. Now set
\begin{equation}
\widetilde{X}_{t} = \sum_{i=1}^{N} x_{i}(t) \mbox{\textbf{1}}_{C_{i}(\chi_{i})}(W), \ \ N \geq 1.
\label{BrownDiffQuant}
\end{equation}

\noindent The process $\widetilde{X}_{t}$ is a non-Voronoi quantization but it is easily computable once the above integral equations are solved since the weights $\mathbb{P}\left(\widehat{W} = \chi_{\underline{i}} \right)$ are known. However the ODE \eqref{eq13} has no explicit solution in general. Then, for numerical implementation purposes, we use discretization schemes like Runge-Kutta one to estimate these quantizers. One shows that the quantized process $\widetilde{X}$ converges toward the process $X$ with respect to the quadratic norm and the rate of convergence is given by the following result.

\begin{prop}[See \cite{Luschgy2006}]
Assume that $b$ is differentiable, $\sigma$ is positive twice differentiable and that $b'-b\frac{\sigma'}{\sigma}-\frac{1}{2}\sigma\sigma''$ is bounded. Then
$$
\left\|X-\widetilde{X}\right\|_{2} =  O\left((\log N)^{-\frac{1}{2}}\right).
$$
\end{prop}

\subsection{Quadrature formulae for numerical integration}
We conclude this first section on optimal (quadratic) quantization by illustrating how to use it for numerical integration of functions defined on the Hilbert space $E$. We provide some quadrature formulae using the above quantization errors. We refer to \cite{Pag`es2005} for the proofs. The main idea is that we know that $\widehat{X}$ is close to $X$ in distribution and if one has a numerical access to the $N$-quantizer $x$ with the associated weights sequence $\left(\mathbb{P}\left(X\in C_{i}(x)\right)\right)_{1 \leq i \leq N}$ of the quantization $\widehat{X}$ then for every Borel functional $F:H\rightarrow \mathbb{R}$, the computation of the expectation
$$
\mathbb{E}\left[F(\widehat{X})\right] = \sum_{i=1}^{N} F(x_{i}) \mathbb{P}\left(X\in C_{i}(x)\right)
$$ 

\noindent is straightforward. The proposition below gives some error bounds for $\mathbb{E}\left[F(X)\right] - \mathbb{E}[F(\widehat{X})]$ based on $L^{p}$-quantization error ($p=2$ or 4) of  $\left\|X-\hat{X}\right\|_{p}$.

\medskip

Let $x$ be a stationary quantizer for $X$ with $\widehat{X}$ its associated Voronoi quantization and $F:E \rightarrow \mathbb{R}$ be a Borel functional defined on $E$.
\begin{itemize}

\item[(i)]Inequality for convex functionals:  If $F$ is convex then
$$
\mathbb{E}\left[F(\widehat{X})\right] \leq \mathbb{E}\left[F(X)\right].
$$

\item[(ii)]Lipschitz functionals: 
\begin{itemize}

\item If $F$ is Lipschitz continuous then
$$
\left|\mathbb{E}\left[F(X)\right] - \mathbb{E}\left[F(\widehat{X})\right]\right| \leq [F]_{Lip} \left\|X-\widehat{X}\right\|_{2}.
$$

\item Let $\theta: E\rightarrow \mathbb{R}_{+}$ be a nonnegative convex function such that $\theta(X) \in L^{2}(\mathbb{P})$. If $F$ is locally Lipschitz with at most $\theta$-growth, $i.e.$ $\left|F(x)-F(y)\right| \leq [F]_{Lip} |x-y| \left(\theta(x)+\theta(y)\right)$ then $F(X)\in L^{1}(\mathbb{P})$ and
$$
\left|\mathbb{E}\left[F(X)\right] - \mathbb{E}\left[F(\widehat{X})\right]\right| \leq 2 [F]_{Lip} \left\|X-\widehat{X}\right\|_{2}\left\|\theta(X)\right\|_{2}.
$$
\end{itemize}

\item[(iii)]Differentiable functionals: If $F$ is differentiable on $E$ with an $\alpha$-H\"older differential D$F$ ($\alpha \in [0,1]$), then
$$
\left|\mathbb{E}\left[F(X)\right] - \mathbb{E}\left[F(\widehat{X})\right]\right| \leq [DF]_{\alpha} \left\|X-\widehat{X}\right\|_{2}^{1+\alpha}.
$$
\end{itemize}

Other quadrature formulae can be derived based on regularity assumptions on $F$ (for more details we refer to \cite{Pag`es2005}).
\section{Quantized importance sampling algorithm}

\subsection{The finite-dimensional setting}
In order to derive the existence of a unique minimum for the function $Q$, we make the following assumption:
\begin{assump}
$$
p \mbox{ is strictly log-concave and} \ \ \lim_{|x| \rightarrow +\infty }p(x)=0.
$$
\end{assump}

\noindent Moreover, the differentiation of the quadratic norm $Q$ (with respect to $\theta$) defined by \eqref{eq2} is required further on and we need for this purpose the following assumptions on the probability density function $p$:
\begin{assump}
The density function $p$ is twice differentiable and satisfies for some $\alpha\in (0,1]$
\begin{itemize}

\item[(i)] $\left|\frac{D p}{p}\right| = O(|x|^{\alpha})$ as $x \rightarrow +\infty $.

\item[(ii)] $\exists \ C>0$ such that $\forall x\in \mathbb{R}^{d}$, $\frac{1}{p(x)}\left|D^{2}p\right|(x) \leq C\left(|x|^{2\alpha}+\frac{1}{|x|^{1-\alpha}}\right)$, where $D^{2}p$ is the Hessian of $p$. 
\end{itemize}
\end{assump}

\begin{prop}
Suppose that Assumptions 1 and 2 are satisfied and the function $F$ satisfies
$$
\forall \ \theta \in \mathbb{R}^{d}, \ \ \mathbb{E}\left[F^{2}(X) \frac{p(X)}{p(X-\theta)}\right] < + \infty \ \ \mbox{ and } \ \ \forall \ M>0, \ \mathbb{E}\left[F^{2}(X)|X|^{\alpha} e^{M|X|^{\alpha}}\right]<+\infty.
$$

\noindent Then, the function $Q$ defined by \eqref{eq2} is finite, strictly convex, differentiable on $\mathbb{R}^{d}$, goes to infinity as $|\theta|$ goes to infinity. As a consequence, the function $Q$ admits a unique global minimum $\theta^{*}$ satisfying
$$
\Argmin Q = \left\{\theta \in \mathbb{R}^{d} \ | \ D Q(\theta) = 0 \right\} = \left\{\theta^{*}\right\}
$$

\noindent and the gradient is given by 
\begin{equation}
D Q(\theta) = \mathbb{E}\left[F^{2}(X)\frac{p(X)}{p^{2}(X-\theta)}D
p(X-\theta)\right].
\label{gradQ}
\end{equation}

\noindent Moreover, if $F$ satisfies 
\begin{equation}
\forall \ M>0, \ \mathbb{E}\left[F^{2}(X)|X|^{2\alpha} e^{M|X|^{\alpha}}\right]<+\infty  
\label{HessHyp1} 
\end{equation}

\noindent and, 
\begin{equation}
\forall \ M>0, \ \sup_{|\theta| \leq M}\mathbb{E}\left[F^{2p}(X+\theta)e^{M|X|^{\alpha}}\frac{1}{|X|^{p(1-\alpha)}}\right]<+\infty, \ \ \mbox{ for some } p>1,
\label{HessHyp2}
\end{equation}

\noindent then, $Q$ is twice differentiable and its hessian is given by
\begin{equation}
D^{2}Q(\theta) = \mathbb{E}\left[F^{2}(X) \frac{p(X)}{p(X-\theta)}\left(2\frac{D p(X-\theta) D p(X-\theta)^{T}}{p^{2}(X-\theta)}-\frac{D^{2}p(X-\theta)}{p(X-\theta)}\right)\right].
\label{HessQ}
\end{equation}
\end{prop}

\begin{proof}[\textnormal{\textbf{Proof.}}]
For every $x\in \mathbb{R}^{d}$, $\theta \mapsto \log p(x-\theta)$ is strictly concave so that $-\log p(x-.)=\log(\frac{1}{p(x-.)})$ is strictly convex, hence the function $\theta \mapsto e^{\log(\frac{1}{p(x-\theta)})}=\frac{1}{p(x-\theta)}$ is strictly convex. Combining Fatou's lemma and Assumption 1, one easily obtains that $\lim_{|\theta| \rightarrow +\infty}Q(\theta)=+\infty$.

\medskip

In order to get the formal differentiation representation \eqref{gradQ} we have to check the domination property for $\theta \in B(0,R)$, for every $R>0$. The log-concavity of $p$ implies that for every $x\in \mathbb{R}^{d}$ and $\theta \in B(0,R)$,
$$
\log p(x) \leq \log p(x-\theta) + \frac{\left\langle D p(x-\theta), \theta \right\rangle}{p(x-\theta)}.
$$

\noindent Hence using Assumption 2 (i) yields
\begin{equation}
0 \leq \frac{p(x)}{p(x-\theta)}  \leq e^{\frac{|D p|}{p}(x-\theta)|\theta|} \leq C_{R}e^{C_{R}|x|^{\alpha}}
\label{logconcavineq}
\end{equation}

\noindent so that 
$$
F^{2}(X)\frac{p(X)}{p^{2}(X-\theta)}D
p(X-\theta) \leq C_{R} F^{2}(X)(|X|^{\alpha}+1)e^{C_{R}|X|^{\alpha}} \in L^{1}(\mathbb{P}).
$$

\noindent To justify the formal differentiation of $DQ$ to get \eqref{HessQ} we proceed as follows. Let $x\in \mathbb{R}^{d}$. Using Assumption 2 yields for every $\theta \in B(0,R)$
$$
\frac{p(x)}{p(x-\theta)}\left(\frac{\left|D p(x-\theta)D p(x-\theta)^{T}\right|}{p^{2}(x-\theta)}\right) \leq C_{R}e^{C_{R}|x|^{\alpha}}\left(|x|^{2\alpha}+1\right)
$$

\noindent and,
$$
\frac{p(x)}{p(x-\theta)} \left(\frac{1}{p(x-\theta)}\left|D^{2}p\right|(x-\theta)\right) \leq C_{R}\left( e^{C_{R}|x|^{\alpha}}|x|^{2\alpha}+\frac{e^{C_{R}|x-\theta|^{\alpha}}}{|x-\theta|^{1-\alpha}}\right).
$$

\noindent Consequently, $\forall \ \theta \in B(0,R)$, we have
$$
F^{2}(x) \frac{p(x)}{p(x-\theta)}\left(2\frac{D p(x-\theta) D p(x-\theta)^{T}}{p^{2}(x-\theta)}-\frac{D^{2}p(x-\theta)}{p(x-\theta)}\right) \leq C_{R}(f(x)+ g_{\theta}(x))
$$

\noindent where $f(x) = F^{2}(x) e^{C_{R}|x|^{\alpha}}\left(|x|^{2\alpha}+1\right)$ and $g_{\theta}(x)=F^{2}(x)\frac{1}{|x-\theta|^{1-\alpha}}e^{C_{R}|x-\theta|^{\alpha}}$. Using the assumption \eqref{HessHyp1} implies that $f(X)\in L^{1}(\mathbb{P})$. Moreover, using another change of variable and \eqref{logconcavineq} yields 
$$
\mathbb{E}\left[g_{\theta}^{p}(X)\right]=\mathbb{E}\left[F^{2p}(X+\theta)e^{Cp|X|^{\alpha}}\frac{1}{|X|^{p(1-\alpha)}}\frac{p(X+\theta)}{p(X)}\right] \leq \mathbb{E}\left[F^{2p}(X+\theta)e^{C|X|^{\alpha}}\frac{1}{|X|^{p(1-\alpha)}}\right]
$$

\noindent so that \eqref{HessHyp2} implies for every $R>0$
$$
\sup_{\theta \in B(0,R)}\mathbb{E}\left[g_{\theta}^{p}(X)\right] <+\infty.
$$

\noindent Consequently, the family $(g_{\theta}(X))_{\theta \in B(0,R)}$ is $\mathbb{P}$-uniformly integrable. This provides the expected representation \eqref{HessQ}. 
\end{proof}

\noindent \textbf{Examples of distributions}
\begin{itemize}

\item[$\rhd$]\emph{The Normal distribution}
$$
p(x) = (2\pi)^{-\frac{d}{2}} e^{-|x|^{2}/2}, \ \ \ x\in \mathbb{R}^{d}.
$$

\noindent Assumption 1 is clearly satisfied. Moreover, we have $\frac{Dp(x)}{p(x)}=-x$ and $\frac{D^{2}p(x)}{p(x)}=xx^{T}-I_{d}$, where $I_{d}$ is the identity matrix of size $d$, so that Assumption 2 is satisfied with $\alpha =1$. A carefull reading of the proof of Proposition 2.1 shows that assumption \eqref{HessHyp2} is useless and one only needs \eqref{HessHyp1}. Moreover, we have
\begin{align*}
DQ(\theta) & = \mathbb{E}\left[F^{2}(X) e^{\frac{|\theta|^{2}}{2}-\left\langle \theta,X\right\rangle}(\theta-X)\right] ,\\
D^{2}Q(\theta) & = \mathbb{E}\left[F^{2}(X) e^{\frac{|\theta|^{2}}{2}-\left\langle \theta,X\right\rangle} \left(I_{d}+\left(\theta-X\right)\left(\theta-X\right)^{T}\right)\right]
\end{align*}

\item[$\rhd$]\emph{The logistic distribution}
$$
p(x) = \frac{e^{x}}{\left(e^{x}+1\right)^{2}}, \ \ \ x \in \mathbb{R}.
$$

\noindent Assumption 1 is satisfied. Assumption 2 holds with $\alpha=1$. 
\item[$\rhd$]\emph{The hyper-exponential distributions}
$$
p(x) = C_{d,a,\sigma}e^{-\frac{|x|^{a}}{\sigma^{a}}}P(x), \ \ \ x\in \mathbb{R}, \ \ a\in[1,2]
$$

\noindent where $P$ is a positive polynomial function. 
\end{itemize}

\medskip 
The main idea is that we know that $Q$, $DQ$ and $D^{2}Q$ can be approximated by
\begin{align}
\widehat{Q}_{N}(\theta) & := \mathbb{E}\left[F^{2}(\widehat{X})\frac{p(\widehat{X})}{p(\widehat{X}-\theta)}\right] \label{QQuant} \\
D\widehat{Q}_{N}(\theta)     & := \mathbb{E}\left[F^{2}(\widehat{X})\frac{p(\widehat{X})}{p^{2}(\widehat{X}-\theta)}D
p(\widehat{X}-\theta)\right] \label{JacobQuant}, \\
D^{2}\widehat{Q}_{N}(\theta) & := \mathbb{E}\left[F^{2}(\widehat{X}) \frac{p(\widehat{X})}{p(\widehat{X}-\theta)}\left(2\frac{D p(\widehat{X}-\theta) D p(\widehat{X}-\theta)^{T}}{p^{2}(\widehat{X}-\theta)}-\frac{D^{2}p(\widehat{X}-\theta)}{p(\widehat{X}-\theta)}\right)\right] \label{HessQuant},
\end{align}

\noindent using an $N$-quantizer $x$ with the associated weights sequence $\left(\mathbb{P}\left(X\in C_{i}(x)\right)\right)_{1 \leq i \leq N}$ of the quantization $\widehat{X}$. The computations of \eqref{QQuant}, \eqref{JacobQuant} and \eqref{HessQuant} are straightforward. For $N$ large enough, $F(x_{i})\neq0$ for some $i \in \left\{1, \cdots, N\right\}$, so that $\widehat{Q}_{N}$ is also strictly convex and goes to infinity as $|\theta|$ goes to infinity. Moreover, it is clear that $\widehat{Q}_{N}$ is differentiable. Hence, there exists a unique $\widehat{\theta}^{N} \in \mathbb{R}^{d}$ such that $\Argmin \widehat{Q}_{N} = \left\{\theta \in \mathbb{R}^{d} \ | \  D\widehat{Q}_{N}(\theta) = 0 \right\} = \left\{\widehat{\theta}^{N}\right\}$. Moreover, for $N$ large enough, the Hessian matrix $D^{2}\widehat{Q}_{N}(\theta)$ is symmetric positive definite for every $\theta \in \mathbb{R}^{d}$.

The following proposition describes the asymptotic behavior of $\widehat{\theta}^{N}$ as $N\rightarrow + \infty$. It shows, as expected, that $\widehat{\theta}^{N}\rightarrow \theta^{*}$ as $N\rightarrow + \infty$. First we need the following assumption:

\begin{assump}
The function $F$ is positive, convex on $\mathbb{R}^{d}$ and satisfies
\begin{itemize}

\item[(i)] F is Lipschitz on $\mathbb{R}^{d}$,

\item[(ii)] $\exists \ a>1$, such that $\mathbb{E}\left[F^{4a}(X)\right]<+\infty$ and for all $M>0$, $\mathbb{E}\left[|X|^{4\frac{a}{a-1}}+ e^{M |X|}\right]<+\infty$.

\end{itemize}
\end{assump}

\begin{prop}[Convergence of $(\widehat{\theta}^{N})_{N\geq1}$]
Consider an $L^{2}$-optimal stationary quantizer $x$ of size $N$ with its associated quantization $\widehat{X}$. Assume that the assumptions of Proposition 2.1 and that Assumption 3 are satisfied. Then, we have 
$$
\widehat{\theta}^{N} \rightarrow \theta^{*} \ \ \mbox{ and } \ \ \widehat{Q}_{N}(\widehat{\theta}^{N}) \rightarrow Q(\theta^{*}) \mbox{ as } N\rightarrow +\infty, 
$$

\noindent where $\widehat{\theta}^{N}$ is the unique global minimum of $\widehat{Q}_{N}$ defined by \eqref{QQuant}.
\end{prop}

\begin{proof}[\textnormal{\textbf{Proof}}]
The first step of the proof consists in showing that the function $\widehat{Q}_{N}$ converges locally uniformly to the continuous function $Q$. Let $x\in \mathbb{R}^{d}$ and $x' \in \mathbb{R}^{d}$. We have
\begin{align*}
\left|F^{2}(x)\frac{p(x)}{p(x-\theta)}-F^{2}(x')\frac{p(x')}{p(x'-\theta)}\right|  & = \left|\left(F^{2}(x)-F^{2}(x')\right)\frac{p(x)}{p(x-\theta)}+F^{2}(x')\left(\frac{p(x)}{p(x-\theta)}-\frac{p(x')}{p(x'-\theta)}\right)\right| \\
  & \leq [F]_{Lip} \left|x-x'\right| \left(|F(x)|+|F(x')|\right)\frac{p(x)}{p(x-\theta)} \\ 
  & \hspace*{5cm} + F^{2}(x')\left|\frac{p(x)}{p(x-\theta)}-\frac{p(x')}{p(x'-\theta)}\right|.
\end{align*}

\noindent Using the log-concavity of $p$, Assumption 2 and the inequality $|e^{u}-e^{v}|\leq \left|u-v\right| (e^{u}+e^{v})$ yields, for every $\theta \in B(0,R)$
\begin{align*}
\left|\frac{p(x)}{p(x-\theta)}-\frac{p(x'-\theta)}{p(x'-\theta)}\right| & \leq \left|\log\left(\frac{p(x)}{p(x')}\right)-\log\left(\frac{p(x-\theta)}{p(x'-\theta)}\right)\right| \left(\frac{p(x)}{p(x')}+ \frac{p(x-\theta)}{p(x'-\theta)}\right) \\
& \leq \left(\left|\log\left(\frac{p(x)}{p(x')}\right)\right| + \left|\log\left(\frac{p(x-\theta)}{p(x'-\theta)}\right)\right|\right) \left(\frac{p(x)}{p(x')}+ \frac{p(x-\theta)}{p(x'-\theta)}\right) \\
& \leq C_{R} \left|x-x'\right|\left(1 + |x|^{\alpha} + |x'|^{\alpha} \right)\left(e^{C |x|}+ e^{C |x'|}\right) \\
& \leq C_{R} \left|x-x'\right|\left(1 + |x| + |x'| \right)\left(e^{C |x|}+ e^{C |x'|}\right).
\end{align*}

\noindent Consequently, we have
\begin{align*}
\left|F^{2}(x)\frac{p(x)}{p(x-\theta)}-F^{2}(x')\frac{p(x')}{p(x'-\theta)}\right|  & \leq C |x-x'| \left(\left(|F(x)|+|F(x')|\right)e^{C|x|^{\alpha}} \right. \\
& \left. \hspace*{3cm} + F^{2}(x')\left(1 + |x| + |x'| \right)\left(e^{C |x|}+ e^{C |x'|}\right)\right),
\end{align*}

\noindent hence, Schwarz's and H\"older's inequalities implies for every $\theta \in B(0,R)$
\begin{align*}
\left|Q(\theta)-\widehat{Q}_{N}(\theta)\right|  & \leq  \mathbb{E}\left[\left|F^{2}(X)\frac{p(X)}{p(X-\theta)}-F^{2}(\widehat{X})\frac{p(\widehat{X})}{p(\widehat{X}-\theta)}\right|\right]\\
&  \leq C \left\|X-\widehat{X}\right\|_{2} \left( \left\|\left(F(X)+F(\widehat{X})\right)e^{C|X|^{\alpha}}\right\|_{2} \right. \\
& \left. \hspace*{3cm} + \left\|F^{2}(\widehat{X})\left(1 + |X| + |\widehat{X}| \right)\left(e^{C |X|}+ e^{C |\widehat{X}|}\right)\right\|_{2} \right) \\
& \leq C \left\|X-\widehat{X}\right\|_{2} \left(\left(\left\|F(X)\right\|_{2a}+\left\|F(\widehat{X})\right\|_{2a}\right)\left\|e^{C|X|^{\alpha}}\right\|_{2\frac{a}{a-1}}\right. \\
& \left. \hspace*{0cm} + \left\|F^{2}(\widehat{X})\right\|_{2a}\left(1 + \left\|X\right\|_{4 \frac{a}{a-1}} + \left\|\widehat{X}\right\|_{4\frac{a}{a-1}} \right)\left(\left\|e^{C |X|}\right\|_{4\frac{a}{a-1}}+ \left\|e^{C |\widehat{X}|}\right\|_{4\frac{a}{a-1}}\right) \right) \\
\end{align*}

\noindent Now $F^{2a}$ is convex since F is and $u \mapsto u^{2a}$ is increasing and convex on $\mathbb{R}_{+}$. Consequently, by Jensen's inequality
$$
\mathbb{E}\left[F^{2a}(\widehat{X})\right]= \mathbb{E}\left[F^{2a}\left(\mathbb{E}\left[X \ | \widehat{X}\right]\right)\right] \leq \mathbb{E}\left[\mathbb{E}\left[F^{2a}(X) \ | \widehat{X}\right]\right]= \mathbb{E}\left[F^{2a}(X)\right].
$$

\noindent Using similar arguments, we have: $\left\|F^{2}(\widehat{X})\right\|_{2a} \leq \left\|F^{2}(X)\right\|_{2a}$, $\left\|\widehat{X}\right\|_{4\frac{a}{a-1}} \leq \left\|X\right\|_{4\frac{a}{a-1}}$, $\left\|e^{C |\widehat{X}|}\right\|_{4\frac{a}{a-1}} \leq \left\|e^{C |X|}\right\|_{4\frac{a}{a-1}}$.
Finally, we obtain for some positive constant $C$ independent of $\theta$ and $N$
$$
\left|Q(\theta)-\widehat{Q}_{N}(\theta)\right| \leq C \left\|X-\widehat{X}\right\|_{2}\rightarrow 0, \ \mbox{ as } N\rightarrow+\infty.
$$ 

\noindent So that $\widehat{Q}_{N}$ converges locally uniformly to $Q$. 

\bigskip

Now let $\varepsilon>0$. $Q$ being continuous at $\theta^{*}$ and strictly convex
$$
\eta :=\inf_{\theta : |\theta-\theta^{*}|\geq \varepsilon} Q(\theta) - Q(\theta^{*}) >0.
$$

\noindent The local uniform convergence of $\widehat{Q}_{N}$ to $Q$ ensures that 
$$
\exists \ N_{\eta}, \ \forall N\geq N_{\eta}, \ \forall \theta \in \mathbb{R}^{d} \ \mbox{such that} \ |\theta-\theta^{*}| \leq \epsilon, \mbox{ then } \left|\widehat{Q}_{N}(\theta)-Q(\theta)\right| \leq \eta/3.
$$

\noindent Assume that $\theta$ satisfies $|\theta-\theta^{*}|\geq \varepsilon$. The convexity of $\widehat{Q}_{N}$ implies that
$$
\widehat{Q}_{N}\left(\theta^{*} + \varepsilon \frac{\theta-\theta^{*}}{|\theta-\theta^{*}|}\right) \leq \frac{\varepsilon}{|\theta-\theta^{*}|} \widehat{Q}_{N}(\theta)+ \left(1-\frac{\varepsilon}{|\theta-\theta^{*}|}\right) \widehat{Q}_{N}(\theta^{*}), 
$$

\noindent so that,
\begin{align*}
\widehat{Q}_{N}(\theta)-\widehat{Q}_{N}(\theta^{*}) & \geq \frac{|\theta-\theta^{*}|}{\varepsilon}\left(\widehat{Q}_{N}\left(\theta^{*} + \varepsilon \frac{\theta-\theta^{*}}{|\theta-\theta^{*}|}\right) - \widehat{Q}_{N}\left(\theta^{*}\right)\right) \\
 & \geq \frac{|\theta-\theta^{*}|}{\varepsilon} \left(Q\left(\theta^{*} + \varepsilon \frac{\theta-\theta^{*}}{|\theta-\theta^{*}|}\right) - Q\left(\theta^{*}\right)-2\eta/3\right) \\
 & \geq \eta/3.
\end{align*}

\noindent Since $\widehat{\theta}^{N}$ is the unique global minimum of $\widehat{Q}_{N}$, we have  $\widehat{Q}_{N}(\widehat{\theta}^{N})-\widehat{Q}_{N}(\theta^{*}) \leq 0$. Consequently, $|\widehat{\theta}^{N}-\theta^{*}| < \varepsilon$ for $N\geq N_{\eta}$ and $(\widehat{\theta}^{N})_{N \geq 1}$ converges to $\theta^{*}$. Combining the local uniform convergence of $(\widehat{Q}_{N})_{N\geq1}$ to $Q$ and the continuity of $Q$ at $\theta^{*}$, we obtain that $\widehat{Q}_{N}(\widehat{\theta}^{N}) \rightarrow Q(\theta^{*})$, as $N\rightarrow+\infty$. This concludes the proof.
\end{proof}

A classical method for estimating $\widehat{\theta}^{N}$, $i.e.$ for solving the system of nonlinear equations $D \widehat{Q}_{N}(\theta)=0$ is Newton-Raphson's algorithm:
\begin{equation}
\widehat{\theta}_{n+1} = \widehat{\theta}_{n} - D^{2}\widehat{Q}_{N}(\widehat{\theta}_{n})^{-1} D\widehat{Q}_{N}(\widehat{\theta}_{n}), \ \  n\geq0, \ \ \widehat{\theta}_{0} \mbox{ given}.
\label{algo_newton}
\end{equation}

\noindent Newton-Raphson's algorithm is attractive because it converges rapidly from any sufficiently good initial guess under standard assumptions. Indeed, since $\widehat{\theta}^{N}$ is the unique solution of $D \widehat{Q}_{N}(\theta)=0$, $D\widehat{Q}_{N}$ is continuously differentiable on $\mathbb{R}^{d}$ and $D^{2}\widehat{Q}_{N}(\theta)$ is a symmetric positive-definite matrix for all $\theta \in \mathbb{R}^{d}$, then the sequence $(\widehat{\theta}_{n})_{n\geq0}$ defined by \eqref{algo_newton} is known to converge toward $\widehat{\theta}^{N}$ if $\widehat{\theta}_{0}$ is sufficiently close to $\widehat{\theta}^{N}$.

At this stage, it is natural to characterize the rate of convergence of $(\widehat{\theta}^{N})_{N\geq1}$ to $\theta^{*}$. To this end, first we need to obtain some error bounds for $\left|DQ(\theta)-D\widehat{Q}_{N}(\theta)\right|$, $\theta \in B(0,R)$, for some $R>0$.

\begin{assump} The function $\frac{Dp}{p}$ is $\alpha$-H\"{o}lder.
\end{assump}

\begin{prop} Assume that Assumptions 1, 2, 3 and 4 hold. Let $x$ be an $L^{2}$-optimal stationary quantizer of size $N$ and $\widehat{X}$ the associated Voronoi quantization. Then, for every $R>0$, for every $\theta \in B(0,R)$
$$
\left|DQ(\theta)-D\widehat{Q}_{N}(\theta)\right| = O\left(\left\|X-\widehat{X}\right\|_{2}^{\alpha}\right), \ \ N\rightarrow +\infty.
$$ 
%

\noindent Hence, $D\widehat{Q}_{N}$ converges locally uniformly to $DQ$.

\begin{proof}[\textnormal{\textbf{Proof.}}] Let $x\in \mathbb{R}^{d}$, $x' \in \mathbb{R}^{d}$ and $\theta \in B(0,R)$. We have
\begin{align*}
F^{2}(x)\frac{p(x)Dp(x-\theta)}{p^{2}(x-\theta)} - F^{2}(x')\frac{p(x')Dp(x'-\theta)}{p^{2}(x'-\theta)} & = \left(F^{2}(x)-F^{2}(x')\right) \frac{p(x)}{p^{2}(x-\theta)}Dp(x-\theta) \\
& + F^{2}(x')\left(\frac{p(x)}{p(x-\theta)}-\frac{p(x')}{p(x'-\theta)}\right)\frac{Dp(x-\theta)}{p(x-\theta)} \\
& + F^{2}(x')\frac{p(x')}{p(x'-\theta)} \left(\frac{Dp(x-\theta)}{p(x-\theta)}-\frac{Dp(x'-\theta)}{p(x'-\theta)}\right).
\end{align*}

\noindent First we take care of the first term of the above sum. Using Assumption 2 (i), Assumption 3 (i) and \eqref{logconcavineq}, yield
$$
\left|\left(F^{2}(x)-F^{2}(x')\right) \frac{p(x)}{p^{2}(x-\theta)}Dp(x-\theta)\right| \leq C_{R}\left|x-x'\right|\left(F(x)+F(x')\right)(1+|x|^{\alpha})e^{C |x|}.
$$

\noindent Now, we focus on the second term. Using similar arguments than the ones used in the proof of Proposition 2.2 and Assumption 2 (i) yield
$$
F^{2}(x')\left|\frac{p(x)}{p(x-\theta)}-\frac{p(x')}{p(x'-\theta)}\right| \frac{\left|Dp(x-\theta)\right|}{p(x-\theta)} \leq C_{R} |x-x'| F^{2}(x')\left(1+|x|+|x'|\right)\left(e^{C|x|}+e^{C|x'|}\right).
$$

\noindent Finally, using \eqref{logconcavineq} and Assumption 4 for the last term implies
$$
F^{2}(x')\frac{p(x')}{p(x'-\theta)} \left|\frac{Dp(x-\theta)}{p(x-\theta)}-\frac{Dp(x'-\theta)}{p(x'-\theta)}\right| \leq C_{R} |x-x'|^{\alpha}F^{2}(x')e^{C|x'|}
$$

\noindent Using similar arguments than the ones used in the proof of Proposition 2.2 yields
\begin{align*}
\left|DQ(\theta)-D\widehat{Q}_{N}(\theta)\right| & \leq \mathbb{E}\left[\left|F^{2}(X)\frac{p(X)Dp(X-\theta)}{p^{2}(X-\theta)} - F^{2}(\widehat{X})\frac{p(\widehat{X})Dp(\widehat{X}-\theta)}{p^{2}(\widehat{X}-\theta)}\right|\right] \\
& \leq C_{R} \left(\left\|X-\widehat{X}\right\|_{2} + \left\|X-\widehat{X}\right\|_{2}^{\alpha}\right) \\
& \leq C_{R} \left\|X-\widehat{X}\right\|_{2}^{\alpha} \ \ \mbox{ as } N\rightarrow +\infty.
\end{align*}
\end{proof}
\end{prop}

Now we are in position to characterize the convergence rate of $(\widehat{\theta}^{N})_{N\geq 1}$ toward $\theta^{*}$. The following result shows that as expected the error $|\widehat{\theta}^{N}-\theta^{*}|$ is controlled by the quantization error.

\begin{thm}Assume that Assumptions 1, 2, 3 and 4 hold. For every $N \geq 1$, there exists $\eta_{N}>0$ such that, if $|\widehat{\theta}_{0}-\widehat{\theta}^{N}|\leq \eta_{N}$, then the sequence $(\widehat{\theta}_{n})_{n\geq0}$ defined by \eqref{algo_newton} converges to $\widehat{\theta}^{N}$. Moreover, the convergence rate of $(\widehat{\theta}^{N})_{N\geq 1}$ to $\theta^{*}$ is based on $L^{2}$-quantization error in the sense that
$$
|\widehat{\theta}^{N}-\theta^{*}| = O\left(\left\|X-\widehat{X}\right\|_{2}^{\alpha}\right), \ \ \mbox{ as } \ \ N\rightarrow +\infty.
$$
\end{thm}

\begin{proof}[\textnormal{\textbf{Proof.}}]
We define the norm $|\theta|_{*}:=\left|D^{2}Q(\theta^{*})\theta\right|$, for $\theta \in \mathbb{R}^{d}$. Since $D^{2}Q(\theta^{*})$ is non-singular,
$$
\frac{1}{\kappa}|\theta| \leq |\theta|_{*} \leq \kappa |\theta|, \ \ \mbox{ for } \ \theta \in \mathbb{R}^{d},
$$


\noindent where $\kappa:=\max\left(|D^{2}Q(\theta^{*})|;|D^{2}Q(\theta^{*})^{-1}|\right)$. Choose $\delta$ sufficiently small that $ 2 \kappa \delta \left(1+\kappa \delta \right)\leq \frac{1}{2}$. We can choose $\varepsilon >0$ sufficiently small such that if $|\theta-\theta^{*}|\leq \varepsilon$ then
\begin{align}
 & \left|D^{2}Q(\theta)-D^{2}Q(\theta^{*})\right| \leq \delta, \label{eqdev1}\\
 & \left|D^{2}Q^{-1}(\theta)-D^{2}Q^{-1}(\theta^{*})\right| \leq \delta,  \label{eqdev2}\\
 & \left|DQ(\theta)-DQ(\theta^{*})-D^{2}Q(\theta^{*})(\theta-\theta^{*})\right| \leq \delta \left|\theta-\theta^{*}\right|. \label{eqdev3} 
\end{align}

\medskip
 
\noindent The continuity of $D^{2}Q$ at $\theta^{*}$, the non-singularity of $D^{2}Q(\theta^{*})$ and the continuous differentiability of $DQ$ at $\theta^{*}$ ensures the existence of such an $\varepsilon$. The existence of an $\eta_{N}$ such that the sequence $(\widehat{\theta}_{n})_{n\geq0}$ converges to $\widehat{\theta}^{N}$ if $|\widehat{\theta}_{0}-\widehat{\theta}^{N}|\leq \eta_{N}$ follows using the same arguments applied to $\widehat{Q}_{N}$.

The convergence of $(\widehat{\theta}^{N})_{N\geq1}$ to $\theta^{*}$ ensures the existence of $N_{0} \in \mathbb{N}^{*}$ such that if $N\geq N_{0}$, $|\widehat{\theta}^{N}-\theta^{*}|\leq \varepsilon$. Hence \eqref{eqdev1}, \eqref{eqdev2} and \eqref{eqdev3} hold for $\theta=\widehat{\theta}^{N}$, $N\geq N_{0}$. 
%

Now, the recursive algorithm \eqref{algo_newton} can be written
\begin{align}
\widehat{\theta}_{n+1} & = \widehat{\theta}_{n} - D^{2}Q(\widehat{\theta}_{n})^{-1} DQ(\widehat{\theta}_{n}) + r_{n} + s_{n} \label{algo_newton_v2}.
\end{align}

\noindent with $\forall n \in \mathbb{N}$, $r_{n}= D^{2}Q(\widehat{\theta}_{n})^{-1}\left(DQ(\widehat{\theta}_{n})-D\widehat{Q}(\widehat{\theta}_{n})\right)$ and $s_{n}= \left(D^{2}Q(\widehat{\theta}_{n})^{-1}-D^{2}\widehat{Q}(\widehat{\theta}_{n})^{-1}\right) D\widehat{Q}(\widehat{\theta}_{n})$.
\noindent Using the following equality
\begin{align*}
D^{2}Q(\theta^{*}) & \left(\widehat{\theta}_{n}-\theta^{*}  - D^{2}Q(\widehat{\theta}_{n})^{-1}DQ(\widehat{\theta}_{n})\right)  = \left(I_{d}+D^{2}Q(\theta^{*})\left( D^{2}Q(\widehat{\theta}_{n})^{-1}- D^{2}Q(\theta^{*})^{-1}\right)\right) \\
& \hspace*{2cm} \left(\left(D^{2}Q(\widehat{\theta}_{n})-D^{2}Q(\theta^{*})\right)(\widehat{\theta}_{n}-\theta^{*})+\left(DQ(\widehat{\theta}_{n})-DQ(\theta^{*})-D^{2}Q(\theta^{*})(\widehat{\theta}_{n}-\theta^{*})\right)\right)
\end{align*}

\noindent and taking norms in \eqref{algo_newton_v2} yields
\begin{align}
\left|\widehat{\theta}_{n+1}-\theta^{*}\right|_{*}  & \leq \left(1+\left|D^{2}Q(\theta^{*})\right| \left|D^{2}Q(\widehat{\theta}_{n})^{-1}- D^{2}Q(\theta^{*})^{-1}\right|\right)\left(\left|D^{2}Q(\widehat{\theta}_{n})- D^{2}Q(\theta^{*})\right| |\widehat{\theta}_{n}-\theta^{*}| \right. \nonumber \\
&  \left. \hspace*{1cm} + \left|DQ(\widehat{\theta}_{n})-DQ(\theta^{*})-D^{2}Q(\theta^{*})(\widehat{\theta}_{n}-\theta^{*})\right|\right) + \left|D^{2}Q(\theta^{*}) r_n\right| + \left|D^{2}Q(\theta^{*}) s_n\right| \label{devNewton} 
\end{align}

\noindent Let $N\geq N_{0}$ and $\eta_{N}>0$ such that $|\widehat{\theta}_{0}-\widehat{\theta}^{N}| \leq  \eta_{N}$ so that $(\widehat{\theta}_{n})_{n\geq1}$ converges toward $\widehat{\theta}^{N}$. The continuity of $DQ$, $D\widehat{Q}$, $D^{2}Q$ and $D^{2}\widehat{Q}$ at $\widehat{\theta}^{N}$	and the non singularity of $D^{2}Q(\widehat{\theta}^{N})$ and $D^{2}\widehat{Q}(\widehat{\theta}^{N})$ yield
$$
\left|D^{2}Q(\theta^{*}) r_n\right| \rightarrow  \left|D^{2}Q(\theta^{*}) D^{2}Q(\widehat{\theta}^{N})^{-1}(DQ(\widehat{\theta}^{N})-D\widehat{Q}(\widehat{\theta}^{N}))\right|, \ \mbox{ and } \ \ \left|D^{2}Q(\theta^{*}) s_n\right| \rightarrow 0 \  \mbox{ as } \  n \rightarrow + \infty.
$$

\noindent Letting $n$ goes to infinity in \eqref{devNewton} and using \eqref{eqdev1}, \eqref{eqdev2}, \eqref{eqdev3} for $\theta=\widehat{\theta}^{N}$, $N\geq N_{0}$, implies that
$$
\left|\widehat{\theta}^{N}-\theta^{*}\right|_{*} \leq 2 \kappa \delta \left(1+\kappa \delta\right) \left|\widehat{\theta}^{N}-\theta^{*}\right|_{*} +  \left|D^{2}Q(\theta^{*}) D^{2}Q(\widehat{\theta}^{N})^{-1}\right|\left|DQ(\widehat{\theta}^{N})-D\widehat{Q}(\widehat{\theta}^{N})\right|.
$$

\noindent Moreover, for $N\geq N_{0}$, using \eqref{eqdev2} we have 
$$
\left|D^{2}Q(\theta^{*}) D^{2}Q(\widehat{\theta}^{N})^{-1}\right|= \left|D^{2}Q(\theta^{*}) \left(D^{2}Q(\widehat{\theta}^{N})^{-1}-D^{2}Q(\theta^{*})^{-1}\right) + I_{d}\right| \leq \left(1+ \kappa \delta\right).
$$

\noindent Finally, Proposition 2.3 and the choice of $\delta$ yield 
$$
\left|\widehat{\theta}^{N}-\theta^{*}\right|_{*} \leq C \left\|X-\widehat{X}\right\|_{2}^{\alpha}, \ \ N \geq N_{0},
$$

\noindent so that $\left|\widehat{\theta}^{N}-\theta^{*}\right| = O\left(\left\|X-\widehat{X}\right\|_{2}^{\alpha}\right)$ as $N\rightarrow +\infty$.
%
%
\end{proof}

\begin{rem} $\bullet$ When $X$ is a $d$-dimensional gaussian vector, Assumptions 2 and 4 are satisfied with $\alpha=1$ so that if $x$ is an $L^{2}$-optimal quantizer with its associated Voronoi quantization $\widehat{X}$, Theorem 1.1. implies the following error bound in Theorem 2.4
$$
\left|\widehat{\theta}^{N}-\theta^{*}\right| \leq C N^{-\frac{1}{d}}, \ \ \mbox{ as } \ \ N\rightarrow + \infty.
$$

\noindent In practice, only a rough estimate of the optimal change of measure parameter $\theta^{*}$ is needed. According to our numerical results, optimal quantization grids of size $N\sim 200,500$ (depending on the dimension $d$) are enough. Concerning $\widehat{\theta}^{N}$, it can be computed with a high precision by a few steps (usually less than 10) of the Newton-Raphson's optimization procedure \eqref{algo_newton}. 

$\bullet$ For the sake of simplicity we used the classical Lipschitz continuous assumption on $F$ but other error bounds can be derived by replacing Assumption 3 $(i)$ by other smoothness assumption.
\end{rem}

\subsection{Quantized importance sampling for Brownian diffusions}
In this section, we extend the Newton-Raphson's algorithm to the infinite dimensional setting, $i.e.$ the case of path-dependent diffusion. We will rely on the Girsanov transform to play the role of mean translator. To be more precise, we consider a $d$-dimensional It\^{o} process $X$ solution to the stochastic differential equation (S.D.E.) 
\begin{E1}
\mbox{d}X_{t}=b(t,X^{t})\mbox{ }\mbox{d}t+\sigma(t,X^{t})\mbox{ }\mbox{d}W_{t}, \mbox{  } X_{0}=x \in \mathbb{R}^{d},
\end{E1} 

\noindent $W=(W_{t})_{t\in [0,T]}$ being a $q\mbox{-}$dimensional standard Brownian motion and where $X^{t}:=(X_{t\wedge s})_{s\in[0,T]}$ is the stopped process at time t, $b:[0,T]\times\mathcal{C}([0,T],\mathbb{R}^{d}) \rightarrow \mathbb{R}^{d}$, $\sigma:[0,T] \times \mathcal{C}([0,T],\mathbb{R}^{d})\rightarrow \mathcal{M}(d,q)$ are measurable with respect to the canonical predictable $\sigma$-field on $[0,T]\times \mathcal{C}([0,T],\mathbb{R}^{d})$.

Under the following assumption 
\begin{hyp3}
\left\{
	\begin{array}{l}
	(i) \hspace*{0.3cm} b(.,0) \mbox{ and } \sigma(.,0)\mbox{ are continuous, }\\
	(ii) \hspace*{0.2cm} \forall t \in[0,T], \forall x,y \in \mathcal{C}([0,T],\mathbb{R}^{d}), \mbox{ } |b(t,y)-b(t,x)|+||\sigma(t,y)-\sigma(t,x)|| \leq C_{b,\sigma} ||x-y||_{\infty}.
	\end{array}
	\right.
\end{hyp3} 

\noindent strong existence and uniqueness of solutions for $(E_{b,\sigma,W})$ can be proved (for more details, see \cite{Rogers1986}). We aim at devising a robust and automatic Newton-Raphson's algorithm based on functional quantization inspired from Section 2.1 for the computation of 
$$
\mathbb{E}\left[F(X)\right]
$$

\noindent where $F$ is a Borel functional defined on $\mathcal{C}\left([0,T],\mathbb{R}^{d}\right)$ such that 
\begin{equation}
F(X) \in L^{2}\left(\mathbb{P}\right) \ \ \mbox{ and } \ \ \ \mathbb{P}\left(F^{2}(X)>0\right)>0.
\label{hypclassiquediff}
\end{equation}

\noindent In this functional framework, the invariance by translation of the Lebesgue measure \eqref{eq1} is replaced by Girsanov Theorem. We consider a translation process given by $\theta \in L^{2}_{T,q}:=L^{2}([0,T],\mathbb{R}^{q})$ which is slightly less general than the ones used in \cite{Lemaire2008}. Indeed, they considered translation processes of the form $\Theta(t,X^{t})$ defined for every $\xi \in \mathcal{C}\left([0,T],\mathbb{R}^{d}\right)$ and $\theta \in L^{2}([0,T],\mathbb{R}^{p})$ by
$$
\Theta(t,\xi):= \varphi(t,\xi^{t})\theta_{t}, \ \ \mbox{ where } \ \ \ \varphi:[0,T] \times \mathcal{C}\left([0,T],\mathbb{R}^{d}\right) \rightarrow \mathcal{M}(q,p),
$$

\noindent is a prespecified bounded Borel function. In what follows, we can easily adapt to this kind of translation processes but for the sake of simplicity, we prefered to focus on this simple case. 

It follows from Girsanov Theorem that for every $\theta \in L^{2}_{T,q}$
$$
\mathbb{E}\left[F(X)\right] = \mathbb{E}\left[F(X^{(\theta)})e^{-\int_{0}^{T}\left\langle \theta_{s},\text{d}W_{s}\right\rangle-\frac{1}{2}\left\|\theta\right\|^{2}_{L^{2}_{T,q}}}\right]
$$

\noindent where $X^{(\theta)}$ denotes the solution to $\left(E_{b+\sigma\theta,\sigma,W}\right)$. Among all these estimators we want to select the one with the lowest quadratic norm so that we want to solve the following minimization problem
$$
\min_{\theta \in L^{2}_{T,q}}Q(\theta)  \ \ \ \mbox{ where } \ \ \ Q(\theta) := \mathbb{E}\left[F^{2}(X^{(\theta)})e^{-2\int_{0}^{T}\left\langle \theta_{s},\text{d}W_{s}\right\rangle-\left\|\theta\right\|^{2}_{L^{2}_{T,q}}}\right].
$$

\noindent Another Girsanov Theorem yields
\begin{equation}
Q(\theta) := \mathbb{E}\left[F^{2}(X)e^{-\int_{0}^{T}\left\langle \theta_{s},\text{d}W_{s}\right\rangle+\frac{1}{2}\left\|\theta\right\|^{2}_{L^{2}_{T,q}}}\right].
\label{exprQfunct}
\end{equation}

In view of numerical implementation of a Newton-Raphson's algorithm to estimate a minimum of $Q$, we are lead to consider a (non trivial) finite dimensional subspace $E$ of $L^{2}_{T,p}$ spanned by an orthonormal basis $\left(e_{1},\cdots,e_{m}\right)$. Like for the finite dimensional framework, our procedure will be based on the representation (as an expectation) of the first differential $DQ$ and the second differential $D^{2}Q$ of $Q$ on $E$ combined with functional quantization of $\left(E_{b,\sigma,W}\right)$.

\begin{prop} Assume that $\mathbb{E}\left[F(X)^{2+\delta}\right]<+\infty$ for some $\delta>0$ as well as assumptions $(\mathcal{H}_{b,\sigma})$ and \eqref{hypclassiquediff} hold. Then the function $Q$ defined by \eqref{exprQfunct} is finite, strictly convex on $L^{2}_{T,q}$ and 
$$
\lim_{\left\|\theta\right\|^{2}_{L^{2}_{T,q}} \rightarrow +\infty} Q(\theta) = \lim_{\left\|\theta\right\|^{2}_{L^{2}_{T,q}} \rightarrow +\infty, \ \theta \in E} Q(\theta) = +\infty.
$$

\noindent Moreover, the function $Q$ is twice differentiable at every $\theta \in L^{2}_{T,q}$ and for every $\psi \in L^{2}_{T,q}$, $\phi \in L^{2}_{T,q}$
\begin{align}
\left\langle DQ(\theta),\psi\right\rangle_{L^{2}_{T,q}} & = \mathbb{E}\left[F^{2}(X)e^{\Phi(\theta)}\left\langle D\Phi(\theta),\psi\right\rangle_{L^{2}_{T,q}}\right] \label{firstderivativefunct} \\
\left(D^{2} Q(\theta)\psi,\phi\right) & = \mathbb{E}\left[F^{2}(X)e^{\Phi(\theta)}\left(\left\langle D\Phi(\theta),\psi\right\rangle_{L^{2}_{T,q}}\left\langle D\Phi(\theta),\phi\right\rangle_{L^{2}_{T,q}} + \left(D^{2}\Phi(\theta)\psi,\phi\right)\right)\right] \label{secondderivativefunct}
\end{align}

\noindent where $\Phi: \theta \mapsto -\int_{0}^{T}\left\langle \theta_{s},\text{d}W_{s}\right\rangle+\frac{1}{2}\left\|\theta\right\|^{2}_{L^{2}_{T,q}}$ is twice differentiable with 
$$
\left\langle D\Phi(\theta),\psi\right\rangle_{L^{2}_{T,q}} = -\int_{0}^{T}\left\langle \psi_{s}, \text{d}W_{s}\right\rangle + \int_{0}^{T}\left\langle \theta_{s},\psi_{s}\right\rangle ds \ \ \ \mbox{ and } \ \ \
\left(D^{2}\Phi(\theta)\psi ,\phi \right) = \int_{0}^{T}\left\langle \phi_{s},\psi_{s}\right\rangle ds.
$$
\end{prop}

\begin{proof}[\textnormal{\textbf{Proof.}}]
Owing to H\"older's inequality of conjugate exponents $(s,t):=(1+\frac{\delta}{2},1+\frac{2}{\delta})$, we have for every $\theta \in L^{2}_{T,q}$
$$
Q(\theta) \leq \mathbb{E}\left[F(X)^{2+\delta}\right]^{\frac{2}{2+\delta}} \mathbb{E}\left[e^{\frac{t}{2}\left\|\theta\right\|^{2}_{L^{2}_{T,q}}-t\int_{0}^{T} \left\langle \theta_{s}, dW_{s}\right\rangle}\right]^{\frac{1}{t}} = \mathbb{E}\left[F(X)^{2+\delta}\right]^{\frac{2}{2+\delta}} e^{(1+\frac{1}{\delta})\left\|\theta\right\|^{2}_{L^{2}_{T,q}}}< +\infty
$$

\noindent since the Dol\'eans exponential $\left(e^{\int_{0}^{t}\left\langle \theta_{s},dW_{s}\right\rangle-\frac{1}{2}\left\|\theta\right\|^{2}_{L^{2}_{T,q}}}\right)_{t\in[0,T]}$ is a true martingale for any $\theta \in L^{2}_{T,q}$. The function $Q$ is strictly convex since the function $\Phi$ is and exp is strictly increasing and strictly convex.

Now, using the trivial equality
$$
e^{-\int_{0}^{T}\left\langle \theta_{s},\text{d}W_{s}\right\rangle+\frac{1}{2}\left\|\theta\right\|^{2}_{L^{2}_{T,q}}} = \left(e^{-\frac{1}{2}\int_{0}^{T}\left\langle \theta_{s},\text{d}W_{s}\right\rangle+\frac{1}{8}\left\|\theta\right\|^{2}_{L^{2}_{T,q}}}\right)^{2} e^{\frac{1}{4}\left\|\theta\right\|^{2}_{L^{2}_{T,q}}}
$$

\noindent and the reverse H\"older's inequality with conjugate exponents $\left(\frac{1}{3},-\frac{1}{2}\right)$ yields
$$
Q(\theta)  \geq \mathbb{E}\left[F^{\frac{2}{3}}(X)e^{\frac{1}{12}\left\|\theta\right\|^{2}_{L^{2}_{T,q}}}\right]^{3} \mathbb{E}\left[e^{\frac{1}{2}\int_{0}^{T}\left\langle \theta_{s},\text{d}W_{s}\right\rangle-\frac{1}{8}\left\|\theta\right\|^{2}_{L^{2}_{T,q}}}\right]^{-2} \geq \mathbb{E}\left[F^{\frac{2}{3}}(X)\right]^{3}e^{\frac{1}{4}\left\|\theta\right\|^{2}_{L^{2}_{T,q}}}
$$

\noindent by the martingale property of the Dol\'{e}ans exponential. Let $\varepsilon>0$ such that $\mathbb{P}\left(F^{2}(X) \geq \varepsilon\right)>0$. We have $Q(\theta)\geq \varepsilon^{\frac{1}{3}} \mathbb{E}\left[\mbox{\textbf{1}}_{F^{2}(X)\geq\varepsilon}\right]e^{\frac{1}{4}\left\|\theta\right\|^{2}_{L^{2}_{T,q}}}=\varepsilon^{\frac{1}{3}}\mathbb{P}\left(F^{2}(X) \geq \varepsilon\right)^{3}e^{\frac{1}{4}\left\|\theta\right\|^{2}_{L^{2}_{T,q}}}$, so that $Q$ goes to infinity as $\left\|\theta\right\|_{L^{2}_{T,q}}$ goes to infinity.

The random functional $\Phi$ from $L^{2}_{T,q}$ into $L^{r}(\mathbb{P})$ (for all $r\geq1$) is twice differentiable since
$$
\forall \ \theta, \ \psi \in L^{2}_{T,q}, \ \ \ \left|\Phi(\theta+\psi)-\Phi(\theta)-\left\langle D\Phi(\theta),\psi\right\rangle_{L^{2}_{T,q}}\right| \leq \frac{1}{2} \left\|\psi\right\|^{2}_{L^{2}_{T,q}}
$$

\noindent and,
$$
\forall \ \theta, \ \psi, \ \phi \in L^{2}_{T,q}, \ \ \ \left\langle D\Phi(\theta+\phi)-D\Phi(\theta),\psi\right\rangle_{L^{2}_{T,q}} = \left\langle \phi,\psi\right\rangle_{L^{2}_{T,q}} = \left(D^{2}\Phi(\theta)\psi,\phi\right)
$$

\noindent where $\psi \mapsto \left\langle D\Phi(\theta),\psi\right\rangle_{L^{2}_{T,q}}$ is a bounded linear operator from $L^{2}_{T,q}$ into $L^{r}(\mathbb{P})$. Using the Cauchy-Schwarz and the Burkholder-Davis-Gundy inequalities, its operator norm satisfies $|||D\phi(\theta)|||_{L^{2}_{T,q},L^{r}(\mathbb{P})} \leq C_{r}(1+\left\|\theta\right\|^{\frac{1}{2}}_{L^{2}_{T,q}})$. Consequently, by the composition rule, we derive that $\Psi: \theta \mapsto e^{\Phi(\theta)}$ is twice diffentiable from $L^{2}_{T,q}$ into $L^{r}(\mathbb{P})$ for any $r\geq1$. Moreover, for every $\psi, \phi \in L^{2}_{T,q}$, we have $\left\langle D \Psi(\theta), \psi \right\rangle_{L^{2}_{T,q}} = e^{\Phi(\theta)} \left\langle D \Phi(\theta), \psi\right\rangle_{L^{2}_{T,q}}$, and $\left(D^{2}\Psi(\theta) \psi, \phi\right)= e^{\Phi(\theta)} \left(\left\langle D \Phi(\theta), \psi\right\rangle_{L^{2}_{T,q}} \left\langle D \Phi(\theta), \phi\right\rangle_{L^{2}_{T,q}} + \left(D^{2}\Phi(\theta) \psi, \phi\right)\right)$.

We conclude that $\theta \mapsto Q(\theta)=\mathbb{E}\left[F^{2}(X) e^{\Phi(\theta)}\right]$ is twice differentiable with first and second differentials characterized by \eqref{firstderivativefunct} and \eqref{secondderivativefunct}.
\end{proof}

The target of the stochastic algorithm investigated in \cite{Lemaire2008} is the minimum of the restriction of Q on $E$, $\theta^{*}_{E}$ which satisfies $DQ_{|E}(\theta^{*}_{E}) \equiv 0 $. Like for the static framework, in order to approximate $DQ_{|E}(\theta)$ and $D^{2}Q_{|E}(\theta)$ (which can be seen respectively as the $m$-tuple $\left(\left\langle DQ(\theta),e_{i}\right\rangle_{L^{2}_{T,q}}\right)_{1 \leq i \leq m}$ and as an $m\times m$ symetric positive definite matrix $D^{2}Q_{|E}(\theta)=\left(D^{2} Q(\theta)e_{i},e_{j}\right)_{1\leq i \leq m, 1\leq j \leq m}$) for every $\theta \in E$, we consider an (non-Voronoi) $N$-functional quantization $\widetilde{X}$ of $(E_{b,\sigma,W})$ given by \eqref{BrownDiffQuant}. Hence, for every $\theta \in E$, we approximate $Q_{E}(\theta)$, $DQ_{E}(\theta)$ and $D^{2}Q_{E}$ by respectively $\widetilde{Q}_{N}(\theta)$, $D \widetilde{Q}_{N}(\theta)$ and $D^{2} \widetilde{Q}_{N}(\theta)$ defined by
\begin{align*}
\widetilde{Q}_{N}(\theta) & = \mathbb{E}\left[F^{2}(\widetilde{X})e^{-\int_{0}^{T}\left\langle \theta_{s},\text{d}\widehat{W}_{s}\right\rangle+\frac{1}{2}\left\|\theta\right\|^{2}_{L^{2}_{T,q}}}\right], \\
\left\langle D \widetilde{Q}_{N}(\theta),e_{i}\right\rangle_{L^{2}_{T,q}} & = \mathbb{E}\left[F^{2}(\widetilde{X})e^{\widehat{\Phi}(\theta)}\left\langle D\widehat{\Phi}(\theta),e_{i}\right\rangle_{L^{2}_{T,q}}\right], \\
\left(D^{2} \widetilde{Q}_{N}(\theta)e_{i},e_{j}\right) & = \mathbb{E}\left[F^{2}(\widetilde{X})e^{\widehat{\Phi}(\theta)}\left(\left\langle D\widehat{\Phi}(\theta),e_{i}\right\rangle_{L^{2}_{T,q}}\left\langle D\widehat{\Phi}(\theta),e_{j}\right\rangle_{L^{2}_{T,q}} + \left(D^{2}\widehat{\Phi}(\theta)e_{i},e_{j}\right)\right)\right]
\end{align*}

\noindent where $\widehat{\Phi}(\theta) = -\int_{0}^{T}\left\langle \theta_{s},\text{d}\widehat{W}_{s}\right\rangle+\frac{1}{2}\left\|\theta\right\|^{2}_{L^{2}_{T,q}}$, $\left\langle D\widehat{\Phi}(\theta),e_{i}\right\rangle_{L^{2}_{T,q}} = -\int_{0}^{T}\left\langle e_{i}(s), \text{d}\widehat{W}_{s}\right\rangle + \int_{0}^{T}\left\langle \theta_{s},e_{i}(s)\right\rangle ds$ and $\left(D^{2}\widehat{\Phi}(\theta)e_{i} ,e_{j}\right) = \left\langle e_{i},e_{j}\right\rangle_{L^{2}_{T,q}}$, $i=1,\cdots, m$, $j=1,\cdots,m$.

Hence, we compute the minimum $\widetilde{\theta}^{N}$ of $\widetilde{Q}_{N}$ by devising the following Newton-Raphson algorithm
\begin{equation}
\widetilde{\theta}_{n+1} = \widetilde{\theta}_{n} - H(\widetilde{\theta}_{n})^{-1} J(\widetilde{\theta}_{n}), \ \ n\geq0,  \ \ \widetilde{\theta}_{0}\in E \ \mbox{ given}.
\label{AlgoNewtonFunct}
\end{equation}

\noindent where for every $\theta \in E$, $H(\theta)_{i,j} = \left(D^{2} \widetilde{Q}_{N}(\theta)e_{i},e_{j}\right)$ and $J(\theta)_{i} = \left\langle D \widetilde{Q}_{N}(\theta),e_{i}\right\rangle_{L^{2}_{T,q}}$, $i=1,\cdots,m$, $j=1,\cdots,m$.
\medskip

\begin{rem}
Like for the static framework (see Proposition 2.2 and Theorem 2.4), the convergence of $\left(\widetilde{\theta}^{N}\right)_{N\geq1}$ toward $\theta^{*}_{E}$ and its convergence rate can be established under assumptions similar to the finite dimensional framework.
\end{rem}
\section{Numerical illustrations}
In the following of this section, we illustrate the performance of these generic variance reduction algorithms both in the finite dimensional framework and in the diffusion framework. The numerical simulations are done in Scilab 5.3.
\subsection{Finite dimensional setting}

\textbf{Basket options}:
We consider basket options with payoffs given by $\left(\sum_{i}^{d} w_{i} S^{i}_{T}-K\right)_{+}$ where $(w_{i},\cdots,w_{d})$ is the vector of weights, $K$ denotes the strike, $T$ is the maturity and $S^{i}_{T}$ is the price at maturity of the ith asset. We assume that each of the $d$ assets under the risk-neutral measure has a price given by a Black-Scholes model driven by the vector of independent Brownian motions $W=\left(W^{1}_{t},\cdots, W^{d}_{t}, \ t\geq 0\right)$, 
$$
S^{i}_{t} = S^{i}_{0}e^{(r-\frac{(\sigma^{i})^{2}}{2})t+\sigma^{i}W^{i}_{t}}\stackrel{d}{=} S^{i}_{0}e^{(r-\frac{(\sigma^{i})^{2}}{2})t+\sigma^{i}\sqrt{t}Z^{i}}\ \ \ S_{0} = \left(S^{1}_{0},\cdots,S^{d}_{0}\right),
$$

\noindent where $Z=(Z^{1},\cdots,Z^{d})$ is a Gaussian vector of size $d$. We price this basket option with different values of the number of assets $d$ and the strike $K$. The quantization grids have the same size $N=200$ and the number of Monte-Carlo simulations $n$ is 100,000 in every case. Note that for each value of $d$ and each value of the strike $K$, the prices are computed using the same pseudo-random number generator initialized with the same \emph{seed}. 

The numerical results are reported in Table \ref{tab:resultbaskopt}. In this table, the first two columns correspond to the different values of the dimension $d$ and the strike $K$. The third and fourth columns correspond to the crude Monte-Carlo estimator and its associated variance. The fifth and sixth columns refer to the Monte-Carlo estimator and its variance using the optimal change of measure $\widehat{\theta}^{N}$ computed with our Newton-Raphson's algorithm. 

\begin{table}[ht]
\footnotesize
\caption{\footnotesize Basket option in dimension $d=2,\cdots,6$ with $r=0.05$, $T=1$, $S_{0}^{i}=50$, $\sigma^{i} = 0.3$, $w_{i}=1/d$, $i=1,\cdots,d$, $N=200$ $n=100,000$}
\centering 
\begin{tabular}{c c c c c c} 
\hline
\textbf{d} & \textbf{K} & \textbf{Price MC} & \textbf{Variance MC} & \textbf{Price QIS} & \textbf{Variance QIS}\\ [0.5ex]
\hline 
2 & 50 & 5.475 & 62.26 & 5.490 & 7.86	\\		
  & 55 & 3.294 & 40.54 & 3.309 & 4.33 \\
  & 60 & 1.873 & 24.03 & 1.885 & 2.03 \\
3 & 50 & 4.751 & 42.01 & 4.760 & 5.88 \\
  & 55 & 2.523 & 24.26 & 2.545 & 2.82 \\
  & 60 & 1.237 & 12.38 & 1.221 & 1.01 \\
4 & 50 & 4.333 & 32.03 & 4.343 & 4.64 \\
	& 55 & 2.086 & 16.93 & 2.089 & 2.01 \\
	& 60 & 0.882 & 7.29  & 0.868 & 0.58 \\
5 & 50 & 4.061 & 26.05 & 4.057 & 3.99 \\
  & 55 & 1.781 & 12.77 & 1.777 & 1.51 \\
  & 60 & 0.642 & 4.63  & 0.647 & 0.35 \\
6 & 50 & 3.843 & 22.10 & 3.830 & 3.50 \\
  & 55 & 1.566 & 10.02 & 1.553 & 1.19\\
  & 60 & 0.506 & 3.25  & 0.494 & 0.22 \\[1ex] 
\hline 
\end{tabular}
\label{tab:resultbaskopt}
\end{table}

\normalsize We can see in this example that our Quantization based Importance Sampling algorithm does reduce the variance by a factor varying from 6 up to 15. Note that it does not require any Monte-Carlo simulations to compute the optimal change of measure and unlike most adaptive importance sampling algorithm it does not need specific parameter tuning. One does not have to set up complicated adjustments when using it, it is fully generic and automatic. Hence, it is a very interesting variance reduction procedure to be used in an industrial way. 

\medskip 

\textbf{Spark spread option}: We consider now an exchange option between gas and electricity (called spark spread) with payoff given by $\left(S^{e}_{T}-h_{R}S^{g}_{T}-C\right)_{+}$ where $S^{e}_{T}$ and $S^{g}_{T}$ denote electricity and gas spot prices at maturity $T$, $h_{R}$ is a heat rate and $C$ is the generation cost. This kind of payoff appears in the pricing of power plant. We assume that the dynamic of electricity and gas spot prices follows the SDE:
$$
dS^{j}_{t} = \theta_j\left(\alpha_j-\log S^{j}_{t}\right)S^{j}_{t} dt + \sigma_{j}dW^{j}_{t}, \ \ j=e,g,
$$ 

\noindent where $W^{e}$ and $W^{g}$ are two independent Brownian motions.
The stochastic processes $X^{j} = \log (S^{j})$, $j=g,e$ are Ornstein-Uhlenbeck processes:
$$
dX^{j}_{t} = \theta_{j}(\mu_{j}-X^{j}_{t})dt+\sigma_{j} dW^{j}_{t}, \ \ \ \mbox{where} \ \ \mu_j = \alpha_{j}-\frac{\sigma^{2}_j}{2\theta_j}.
$$

\noindent Writing spot prices as exponential of a sum of Ornstein-Uhlenbeck processes is a very common way to reproduce the mean reversion behavior of commodity spot prices. This model was first proposed by Schwartz in \cite{Schwartz1997}. In this example, the dimension $d$ is equal to 2. The quantization grids have the same size $N=200$ and the number of Monte-Carlo simulations $n$ is 100,000 in every case.
  
The numerical results are summarized in Table \ref{tab:resultsparkspreadopt} where we price the spark spread option for different values of $C$. We see that our Quantization based IS algorithm perfoms well again. In any case, the variance is divided by at least 13. Once again the Newton-Raphson algorithm proposed converges quickly, $i.e.$ five iterations are enough to get a very accurate estimate of $\widehat{\theta}^{N}$.
  
\begin{table}[ht]
\footnotesize
\caption{\footnotesize Spark spread option with $T=0.5$, $S^{e}_{0}=40\ \$ /\textnormal{MWh}$, $S^{g}_{0}=4\ \$ /\textnormal{MMBTU}$ (\textnormal{BTU}: British Thermal Unit), $\sigma_{e} = 0.7$, $\sigma_{g} = 0.35$, $\lambda_e=\lambda_g=0.3$, $\alpha_e=\log(S^{e}_0)$, $\alpha_g=\log(S^{g}_0)$, $h_R=10 \ \textnormal{BTU}/\textnormal{kWh}$, $C=0,3,5,8, 10, 12 \ \$/\textnormal{MWh}$, $N=200$, $n=100,000$.}
\centering 
\begin{tabular}{c c c c c} 
\hline
\textbf{C} & \textbf{Price MC} & \textbf{Variance MC} & \textbf{Price QIS} & \textbf{Variance QIS}\\ [0.5ex]
\hline 
 0  & 7.933 & 221.01 & 7.957 & 16.48	\\
 3  & 6.681 & 189.24 & 6.757 & 13.32	\\
 5  & 6.024 & 176.93 & 6.049 & 11.54	\\
 8  & 5.081 & 153.44 & 5.083 & 9.16	  \\		
 10 & 4.575 & 141.09 & 4.531 & 7.81   \\
 12 & 4.057 & 125.49 & 4.032 & 6.61   \\[1ex] 
\hline 
\end{tabular}
\label{tab:resultsparkspreadopt}
\end{table}

\subsection{Infinite dimensional setting}
We consider three different basis of $L^{2}([0,1],\mathbb{R})$ 
\begin{itemize}

\item a polynomial basis composed of the shifted Legendre polynomials $(\tilde{P}_{n})_{n\geq0}$ defined by 
\begin{shleg} 
\forall n\geq0, \ \forall t \in [0,1], \ \ \tilde{P}_{n}(t)=P_{n}(2t-1) \ \mbox{ where } \ P_{n}(t)=\frac{1}{2^{n} n!}\frac{\mbox{d}^{n}}{\mbox{d}t^{n}}\left(\left(t^{2}-1\right)\right). 
\end{shleg}

\item the Karhunen-Loeve basis defined by \begin{kl}\forall n\geq0, \ \forall t \in [0,1], \ \ e_{n}(t)=\sqrt{2} \ \sin\left(\left(n+\frac{1}{2}\right) \pi t\right).\end{kl}

\item the Haar basis which is defined by 
\begin{haar} 
\forall n \geq0, \ \forall k=0,..., 2^{n}-1, \ \forall t \in [0,1], \ \psi_{n,k}(t)=2^{\frac{k}{2}}\psi(2^{k}t-n) 
\end{haar} 

\noindent where
\begin{eqnarray*}
\psi(t)=\left\{
\begin{array}{l}
1  \ \ \ \  \mbox{ if }t \in[0,\frac{1}{2}) \\
-1 \ \ \mbox{ if }t \in[\frac{1}{2},1) \\ 
0  \ \ \ \ \mbox{ otherwise. }   
\end{array}
\right.
\end{eqnarray*}

\end{itemize}

\textbf{Asian option}: The considered payoff is an Asian option on a discrete time schedule of observation dates $t_{0}<\cdots < t_{p-1}=T$ with payoff $\left(\frac{1}{p}\sum_{k=0}^{p-1}S_{t_{k}} - K\right)_{+}$.

\medskip 

\noindent \emph{Black-Scholes Model:}  First we consider that $S$ follows the classical Black-Scholes model with interest rate $r=4\%$ and volatility $\sigma=50\%$. The strike of the option is set at $K=115$, the maturity at $T=1$, $p=100$ observation dates and the initial price $S_{0}=100$. For the different basis mentioned above and different values of $m$ (2, 4 and 8), the results of our algorithm are summarized in Table \ref{tab:resultAsianBS}.

Note that since for all $t \in [0,T]$, $S_{t} = e^{(r-\frac{\sigma^{2}}{2})t + \sigma W_{t}}$, for numerics we consider the non-Voronoi $\widehat{S}_{t}$ defined by 
\begin{equation}
\widehat{S}_{t} := e^{(r-\frac{\sigma^{2}}{2})t + \sigma \widehat{W}_{t}}= \sum_{i=1}^{N} e^{(r-\frac{\sigma^{2}}{2})t + \sigma \chi_{i}(t)} \mbox{\textbf{1}}_{C_{i}(\chi_{i})}(W), \ \ N \geq 1,
\label{BSNonVoronoiQuant}
\end{equation}

\noindent instead of approximating the solution of the ODE given by \eqref{eq13}.

\begin{figure}[ht]
\begin{center}
   \includegraphics[width=8.0cm,height=6.0cm]{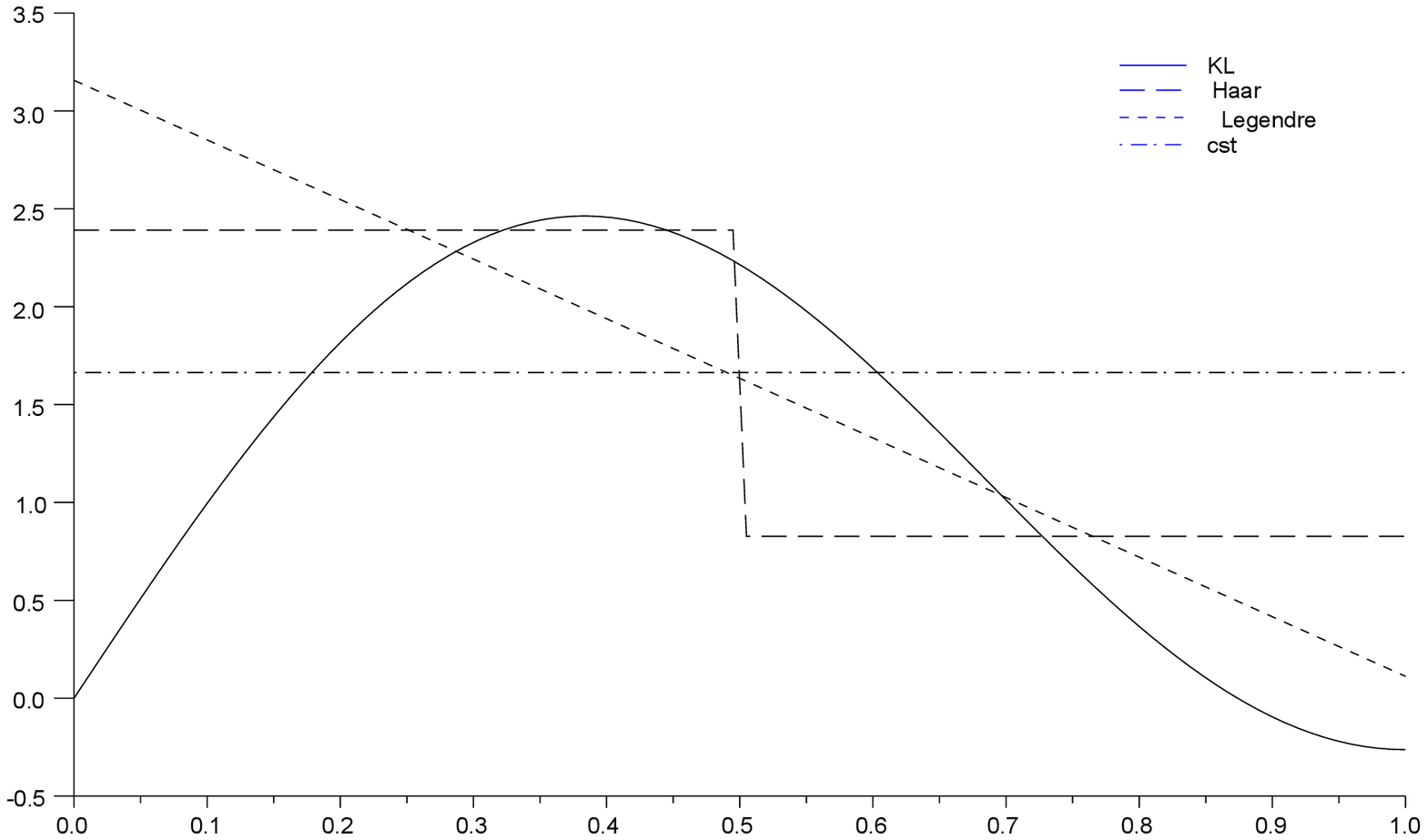}	
   \includegraphics[width=8.0cm,height=6.0cm]{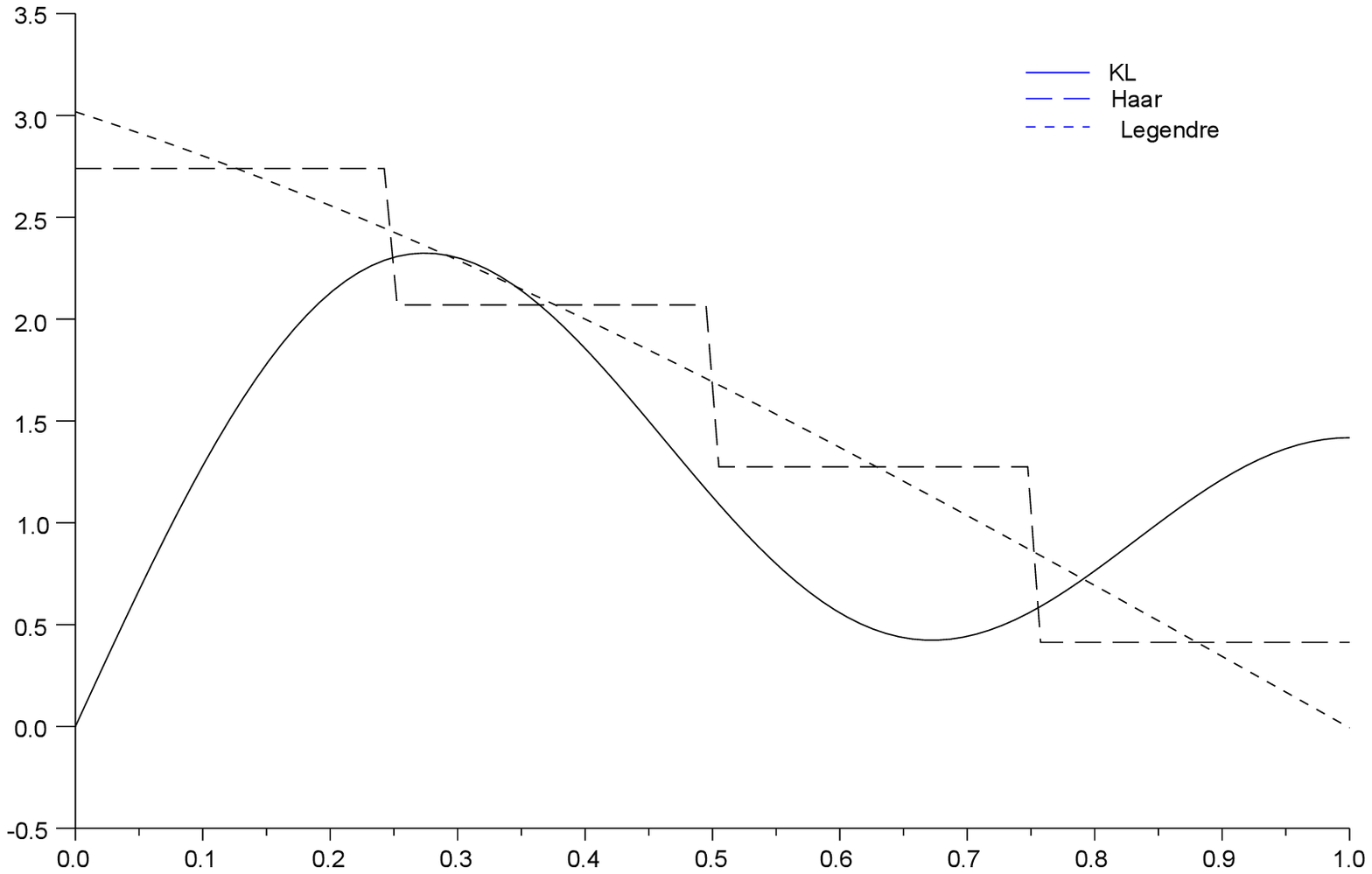}
   \includegraphics[width=8.0cm,height=6.0cm]{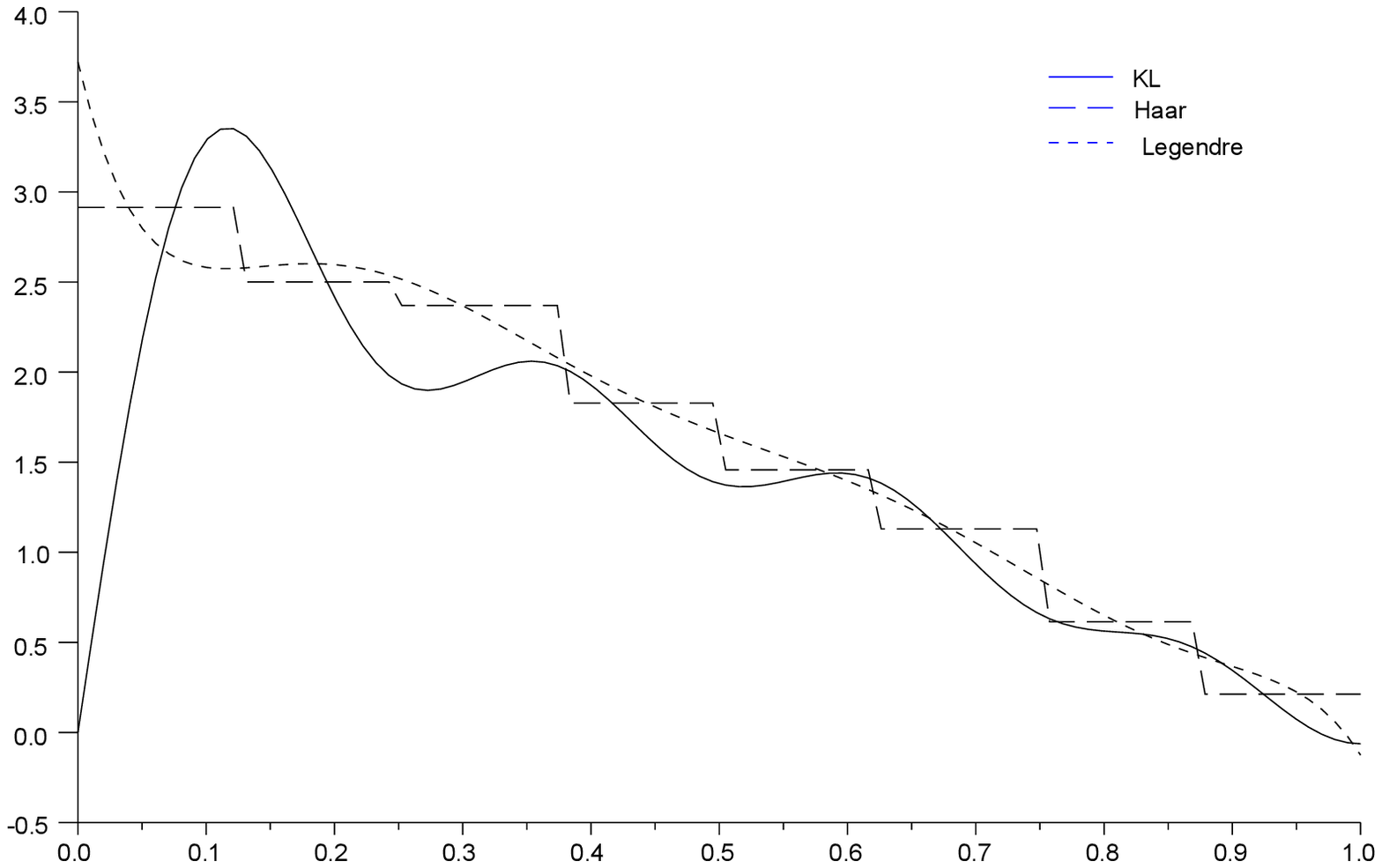}
   \caption{\label{figImpParameterAsianBS}Asian option: Optimal $\theta^{N}$ obtained by our algorithm in the case of the Black-Scholes model for different basis and several values of $m$ ($m=2$ for the left upper curves, $m=4$ for the right upper curves and $m=8$ for the lower curves).}
  \end{center}
\label{fig1}
\end{figure}

We set the optimal product quantizer at level $d_{N}=966$ which corresponds to the optimal decomposition $N_{1}=23$, $N_{2}=7$, $N_{3}=3$, $N_{4}=2$ for the problem \eqref{minprobfunctquant} (see \cite{Luschgy2002} for more details). The number of Monte-Carlo simulations $n$ is 100,000 in every case. Note once again that for each basis and each value of the dimension $m$, the prices and the variances are computed using the same pseudo-random number generator initialized with the same \emph{seed}. In Figure \ref{fig1} are depicted the optimal variance reducer $\theta^{N}$ when the minimization of $\widetilde{Q}_{N}$ is carried out on $E_{m}=\mbox{span}(e_{1},\cdots,e_{m})$ for several values of $m$ in the different basis mentioned above.
\begin{table}[h]
\footnotesize
\caption{\footnotesize Asian option in the Black-Scholes model with $S_0=100$, $K=115$, $T=1$, $p=100$, $\sigma=50\%$, $r=4\%$, $d_{N}=966$, $n=100,000$.}
\centering 
\begin{tabular}{c c c c c c} 
\hline
\textbf{Basis} & \textbf{m} &\textbf{Price MC} & \textbf{Variance MC} & \textbf{Price QIS} & \textbf{Variance QIS}\\ [0.5ex]
\hline 
Constant         & 1  & 7.112  & 293.98 & 7.006 & 79.07	\\
Legendre         & 2  & 7.179  & 296.52 & 7.062 & 22.46	\\
(ShLeg)          & 4  & 7.033  & 290.71 & 7.093 & 22.17	\\
                 & 8  & 7.180  & 300.72 & 7.096 & 22.25	\\
Karhunen-Lo\`eve & 2  & 7.102  & 296.91 & 7.043 & 83.81	\\
(KL)             & 4  & 7.066  & 295.76 & 7.034 & 74.31	\\
                 & 8  & 7.082  & 290.51 & 7.096 & 37.69	\\
Haar						 & 2  & 7.104  & 293.13 & 7.035 & 33.14	\\
(Haar)           & 4  & 7.110  & 297.77 & 7.074 & 25.07	\\
                 & 8  & 7.136  & 299.64 & 7.065 & 23.17	\\[1ex] 
\hline 
\end{tabular}
\label{tab:resultAsianBS}
\end{table}  

\noindent \emph{Local-Volatility Model:} Now, we consider the same payoff function in a local volatility model (inspired by the CEV model) defined by 
\begin{equation}
\mbox{d}X_{t}=rX_{t} \mbox{d}t+ \sigma X_{t} \frac{X_{t}^{\beta}}{\sqrt{1+X_{t}^{2}}} \ \mbox{d}W_{t}, \ X_{0}=x, 
\label{cev}
\end{equation}

\noindent with $r=0.04$, $\sigma=5$, $x=100$ and $\beta=0.5$. 
For numerics, the solution of the ODE given by \eqref{eq13} is approximated by a sixth order Runge-Kutta scheme. The number of Monte-Carlo simulations $n$ is 50,000 in every case. The numerical results are summarized in Table \ref{tab:resultAsianCEV}.
\begin{table}[ht]
\footnotesize
\caption{\footnotesize Asian option in the Local Volatility model with $x=100$, $K=115$, $T=1$, $p=100$, $\sigma=5$, $r=4\%$, $\beta=0.5$, $d_{N}=966$, $n=50,000$.}
\centering 
\begin{tabular}{c c c c c c} 
\hline
\textbf{Basis} & \textbf{m} &\textbf{Price MC} & \textbf{Variance MC} & \textbf{Price QIS} & \textbf{Variance QIS}\\ [0.5ex]
\hline 
Constant         & 1  & 6.635 & 205.69 & 6.681 & 58.63	\\
Legendre         & 2  & 6.646 & 204.84 & 6.619 & 16.94	\\
(ShLeg)          & 4  & 6.593 & 206.52 & 6.661 & 16.84	\\
 								 & 8  & 6.537 & 203.39 & 6.627 & 17.29	\\
Karhunen-Lo\`eve & 2  & 6.562 & 203.73 & 6.620 & 66.55	\\
(KL)             & 4  & 6.627 & 205.24 & 6.578 & 53.25	\\
                 & 8  & 6.700 & 207.31 & 6.637 & 32.61	\\
Haar						 & 2  & 6.583 & 204.65 & 6.651 & 26.06	\\
(Haar)           & 4  & 6.535 & 203.09 & 6.669 & 18.82	\\
                 & 8  & 6.679 & 206.39 & 6.656 & 17.48	\\[1ex] 
\hline 
\end{tabular}
\label{tab:resultAsianCEV}
\end{table}  

\medskip 

\noindent \emph{Schwartz's Model:} Again we consider the same payoff function in the Schwartz model with $r=0.04$, $\sigma=50\%$, $S_{0}=100$, $\alpha=\log(S_{0})$, $\lambda=0.3$. The number of Monte-Carlo simulation $n$ is 100,000 in every case. The numerical results are summarized in Table \ref{tab:resultAsianSchwartz}. Note that since the spot price can be written
$$
S_{t} = e^{\log(S_{0})e^{-\theta t}+\mu(1-e^{-\theta t})+ Y_{t}} \ \mbox{ with } \ \ dY_{t} = -Y_{t} dt +\sigma dW_{t}, \ Y_{0}=0.
$$

Hence, to quantize the diffusion $S$, we just have to obtain a (rate optimal) $N$-product quantizer $y^{N}$ of the centered Ornstein-Uhlenbeck process $Y$. It is given by
$$
y^{N}_{i}(t) = \sigma \sum_{n\geq1} x_{i_{n}} \tilde{c}_{n} \varphi_{n}(t), \ \ \ \ i\in \Pi_{n \geq 1} \left\{1,\cdots,N_{n}\right\},
$$

\noindent where $x_{n}:=\left(x_{1}^{N_{n}},\cdots,x_{N_{n}}^{N_{n}}\right)$ is the unique optimal $N_{n}$-quantizer of the normal distribution, $\tilde{c}_{n} = \frac{T^{2}}{\left(\pi (n-1/2)\right)^{2}+ (\theta T)^{2}}$ ,  $\varphi_{n}(t) = \sqrt{\frac{2}{T}}\left(\frac{\pi}{T}(n-1/2)\sin\left(\pi(n-1/2)\frac{t}{T}\right)+\theta\left(\cos\left(\pi(n-1/2)\frac{t}{T}\right)-e^{-\theta t}\right)\right)$. Consequently we save computation time since we don't need to devise a Runge-Kutta scheme.

\begin{table}[ht]
\footnotesize
\caption{\footnotesize Asian option in Schwartz's model with $S_{0}=100$, $K=115$, $T=1$, $p=100$, $\sigma=50\%$, $r=4\%$, $\alpha=\log(S_{0})$, $\lambda = 0.3$, $d_{N}=966$, $n=100,000$.}
\centering 
\begin{tabular}{c c c c c c} 
\hline
\textbf{Basis} & \textbf{m} &\textbf{Price MC} & \textbf{Variance MC} & \textbf{Price QIS} & \textbf{Variance QIS}\\ [0.5ex]
\hline 
Constant         & 1  & 5.012 & 173.59 & 4.905 & 40.07	\\
Legendre         & 2  & 5.029 & 173.80 & 4.952 & 11.74	\\
(ShLeg)          & 4  & 4.980 & 170.77 & 4.978 & 11.78	\\
 								 & 8  & 5.091 & 180.02 & 4.962 & 12.00	\\
Karhunen-Lo\`eve & 2  & 4.928 & 171.54 & 4.960 & 44.55 	\\
(KL)             & 4  & 4.974 & 171.92 & 4.956 & 39.48	\\
                 & 8  & 4.980 & 171.64 & 4.961 & 21.50 	\\
Haar						 & 2  & 4.999 & 173.55 & 4.944 & 17.47	\\
(Haar)           & 4  & 5.027 & 175.04 & 4.949 & 13.28	\\
                 & 8  & 4.932 & 169.83 & 4.970 & 12.28	\\[1ex] 
\hline 
\end{tabular}
\label{tab:resultAsianSchwartz}
\end{table}  

\medskip

\textbf{Down \& In Call option}: We consider an Down \& In Call option of strike $K$ and barrier $L$. This option is activated when the underlying process $X$ moves down and hits the barrier $L$. The payoff function at maturity T is defined by 
$$
F(X)=(X_{T}-K)_{+}\mbox{\bf 1}_{\left\{ \min_{0 \leq t \leq T} X_{t} \leq L\right\}}
$$

\noindent A standard approach to price the option is to consider the continuous Euler scheme $\bar{X}$ of step $t_{k}=k\frac{T}{M}$ obtained by extrapolation of the Brownian Motion between two instants of discretization. For every $t\in[t_{k},t_{k+1}]$, we can write 
$$
\bar{X}_{t}=\bar{X}_{t_{k}}+b(\bar{X}_{t_{k}})(t-t_{k})+\sigma(\bar{X}_{t_{k}})(W_{t}-W_{t_{k}}), \ \bar{X}_{0}=x_{0} \in \mathbb{R}.
$$ 

\noindent By preconditioning, 
$$
\mathbb{E}\left[(\bar{X}_{T}-K)_{+}\mbox{\bf 1}_{\left\{ \min_{0 \leq t \leq T} \bar{X}_{t} \leq L\right\}}\right]=\mathbb{E}\left[(\bar{X}_{T}-K)_{+}\left(1-\prod_{k=0}^{M-1}p(\bar{X}_{t_{k}},\bar{X}_{t_{k+1}})\right)\right],
$$ 

\noindent where $p(x_{k},x_{k+1})=\mathbb{P}\left(\min_{t_{k}\leq t \leq t_{k+1}}\bar{X}_{t} \geq L\left|\left(\bar{X}_{t_{k}},\bar{X}_{t_{k+1}}\right)=(x_{k},x_{k+1}) \right.\right)$ is the probability of non exit of some brownian bridge. Using the law of the brownian bridge (see for example \cite{Gobet2000}), we can write 
\begin{eqnarray}
p(x_{k},x_{k+1})= 1 - \mathbb{P}\left( \left.\min_{t\in [0,t_{1}]} W_{t} \leq \frac{L-x_{k}}{\sigma(x_{k})} \right| W_{t_{1}} = \frac{x_{k+1}-x_{k}}{\sigma(x_{k})}\right) \\
= \left\{\begin{array}{l}
 1-e^{-\frac{2(L-x_{k})(L-x_{k+1})}{(t_{k+1}-t_{k})\sigma^{2}(x_{k})}} \ \mbox{ if } L\leq \min(x_{k},x_{k+1}), \\
 \\
0 \ \hspace*{3.2cm} \mbox{, otherwise.}
\end{array}
\right.
\label{defprob}
\end{eqnarray}

\noindent Hence, we run our algorithm with this modified payoff function:  $(\bar{X}_{T}-K)_{+}\left(1-\prod_{k=0}^{M-1}p(\bar{X}_{t_{k}},\bar{X}_{t_{k+1}})\right)$. We set the number of steps $M=100$. In the following simulations, we consider the local volatility model \eqref{cev} and the classical Black-Scholes model. The results are summarized in Table \ref{tab:resultDICLocalVol} and Table \ref{tab:resultDICBS}. In Figure 2 are depicted the optimal variance reducer for the local volatility model. Our numerical results illustrate the effectiveness of our Newton-Raphson's IS algorithm. In this example, the computation time needed to achieve a given precision is divided by a factor 8 in comparison with the crude Monte Carlo estimator. 

\begin{table}[ht]
\footnotesize
\caption{\footnotesize Down\&In Call option in Local volatility model with $X_{0}=100$, $K=115$, $L=65$, $T=1$, $\sigma=5$, $r=4\%$, $\beta=0.5$, $d_{N}=966$, $n=50,000$, $M=100$.}
\centering 
\begin{tabular}{c c c c c c} 
\hline
\textbf{Basis} & \textbf{m} &\textbf{Price MC} & \textbf{Variance MC} & \textbf{Price QIS} & \textbf{Variance QIS}\\ [0.5ex]
\hline 
Constant         & 1  & 0.684 & 25.80 & 0.673 & 15.33	\\
Legendre         & 2  & 0.711 & 27.92 & 0.662 & 4.58	\\
(ShLeg)          & 4  & 0.683 & 25.66 & 0.684 & 3.46	\\
								 & 8  & 0.686 & 26.59 & 0.685 & 3.35	\\
Karhunen-Lo\`eve & 2  & 0.680 & 25.38 & 0.696 & 5.22 	\\
(KL)             & 4  & 0.702 & 26.39 & 0.683 & 6.53	\\
  							 & 8  & 0.687 & 26.39 & 0.688 & 5.80	\\
Haar						 & 2  & 0.648 & 24.90 & 0.673 & 8.15	\\
(Haar)           & 4  & 0.671 & 25.15 & 0.692 & 5.46	\\
                 & 8  & 0.709 & 30.17 & 0.700 & 5.29 	\\[1ex] 
\hline 
\end{tabular}
\label{tab:resultDICLocalVol}
\end{table}

\begin{figure}[ht]
\begin{center}
   \includegraphics[width=8.0cm,height=6.0cm]{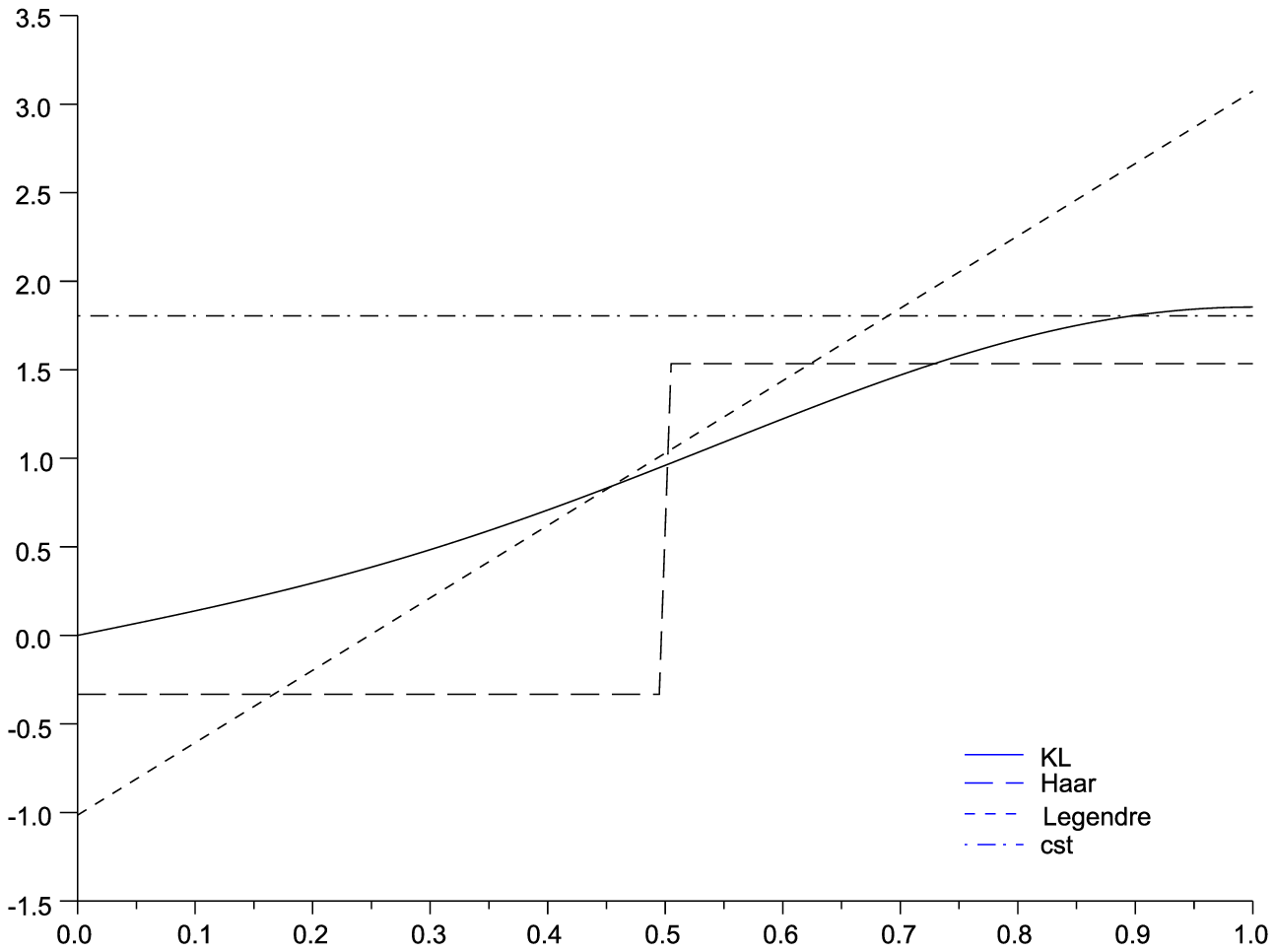}	
   \includegraphics[width=8.0cm,height=6.0cm]{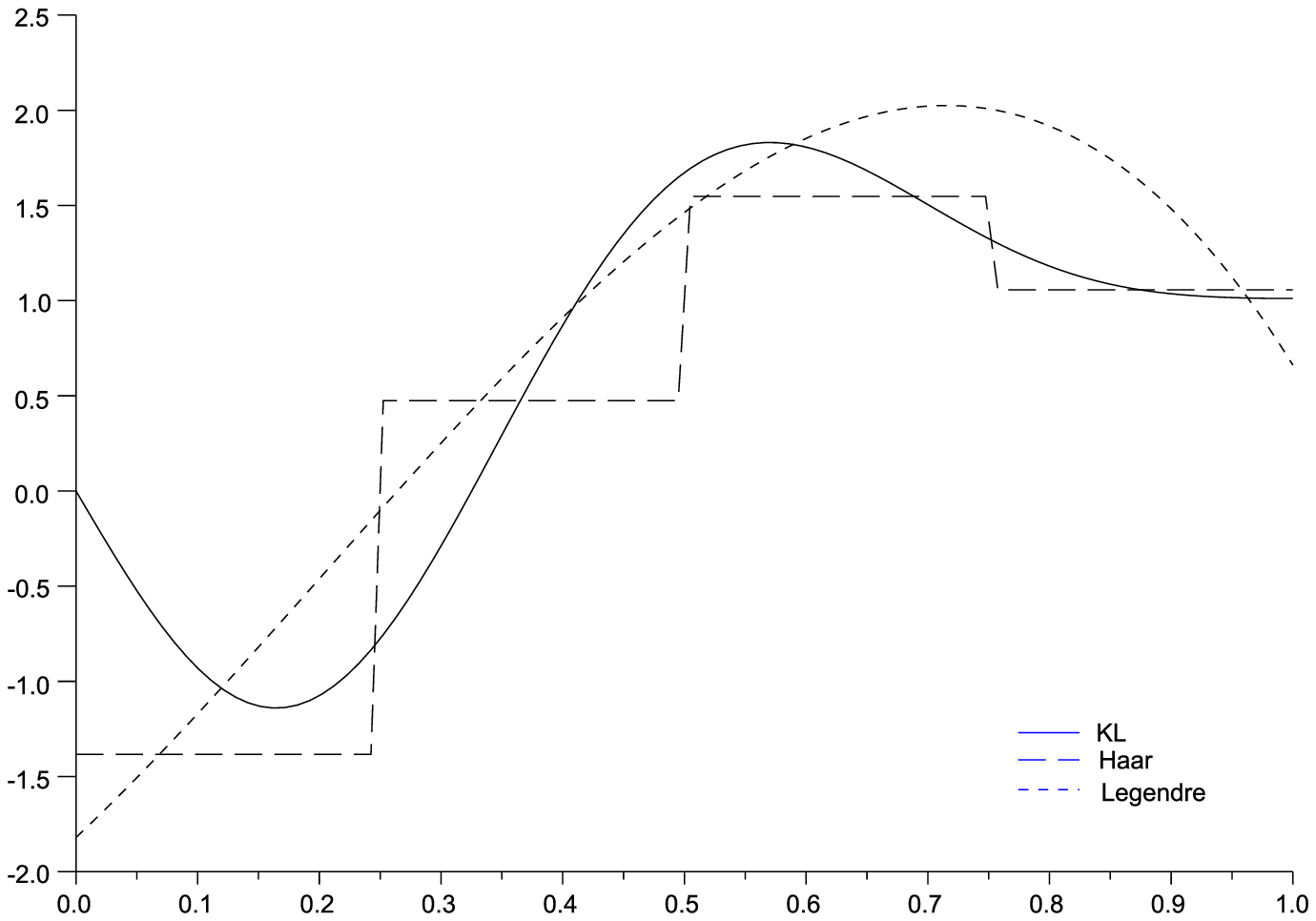}
   \includegraphics[width=8.0cm,height=6.0cm]{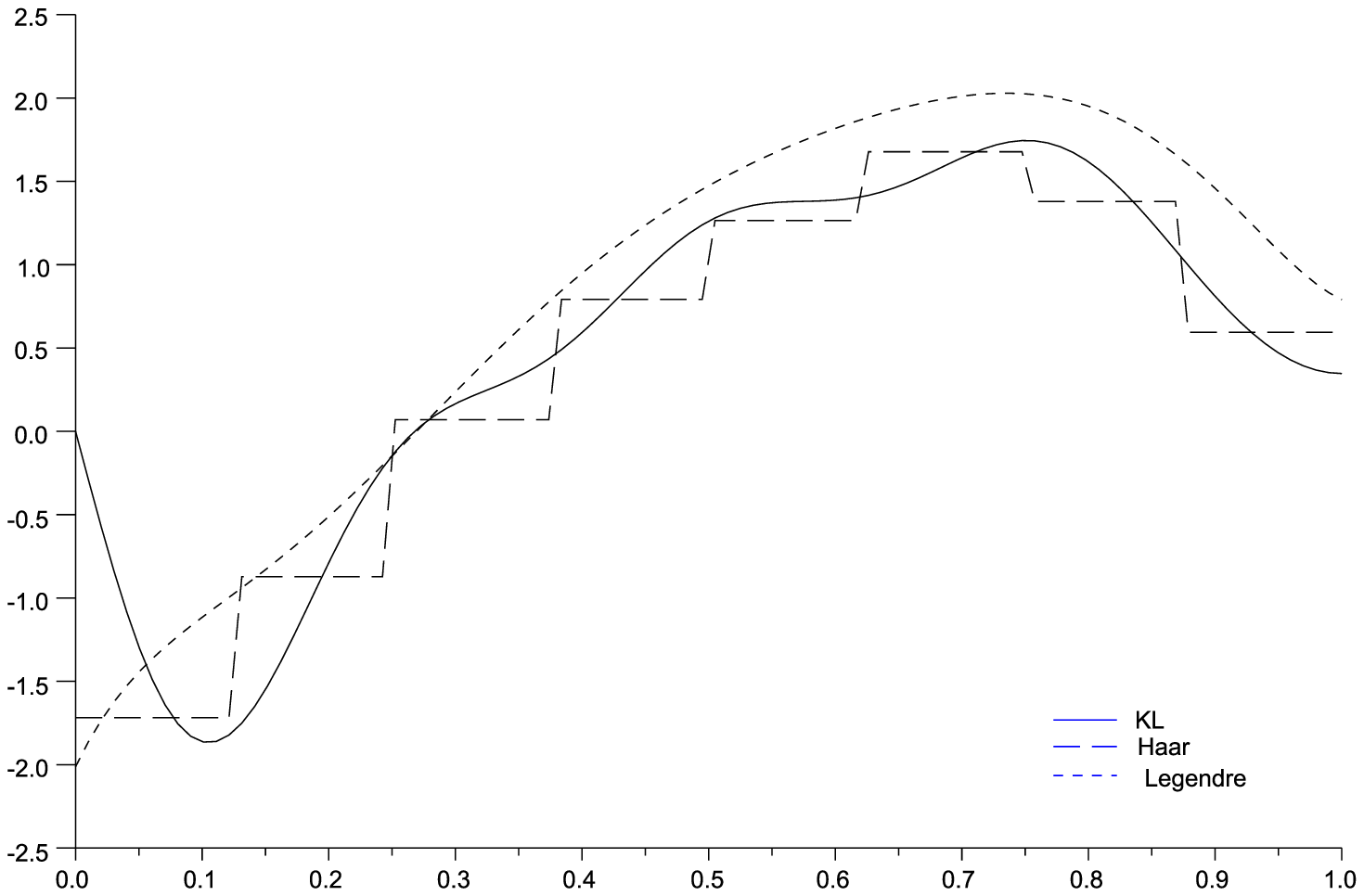}
   \caption{\label{figImpParameterDICLocVol}Down\&In Call option: Optimal $\theta^{N}$ obtained by our algorithm in the case of the Local Volatility model for different basis and several values of $m$ ($m=2$ for the left upper curves, $m=4$ for the right upper curves and $m=8$ for the lower curves).}
  \end{center}
\label{fig2}
\end{figure}

\begin{table}[ht]
\footnotesize
\caption{\footnotesize Down\&In Call option in the Black-Scholes model with $S_{0}=100$, $K=115$, $L=65$, $T=1$, $\sigma=50\%$, $r=4\%$, $d_{N}=966$, $n=100,000$, $M=100$.}
\centering 
\begin{tabular}{c c c c c c} 
\hline
\textbf{Basis} & \textbf{m} &\textbf{Price MC} & \textbf{Variance MC} & \textbf{Price QIS} & \textbf{Variance QIS}\\ [0.5ex]
\hline 
Constant         & 1  & 0.481 & 21.76 & 0.467 & 8.62	\\
Legendre         & 2  & 0.455 & 19.25 & 0.469 & 3.54	\\
(ShLeg)          & 4  & 0.474 & 21.03 & 0.477 & 3.43	\\
								 & 8  & 0.451 & 19.96 & 0.470 & 3.39	\\
Karhunen-Lo\`eve & 2  & 0.466 & 21.36 & 0.459 & 5.51	\\
(KL)             & 4  & 0.470 & 22.37 & 0.471 & 5.78	\\
  							 & 8  & 0.462 & 21.93 & 0.465 & 5.20	\\
Haar						 & 2  & 0.471 & 22.16 & 0.473 & 7.38	\\
(Haar)           & 4  & 0.469 & 22.08 & 0.464 & 5.03	\\
                 & 8  & 0.470 & 21.09 & 0.476 & 5.26	\\[1ex] 
\hline 
\end{tabular}
\label{tab:resultDICBS}
\end{table}

\small
\bibliographystyle{abbrv}
\bibliography{qimportancesamplingbib}

\end{document}